\documentclass{amsart}
\usepackage{amscd, amssymb,epsf,OVmacros}
\begin{document}
\pagenumbering{roman}
\title{Patchworking real algebraic varieties}
\author{Oleg Viro}
\address{Mathematical Department, Uppsala University,
S-751 06 Uppsala, Sweden;\abz
POMI, Fontanka 27, St.~Petersburg, 191011, Russia}
\email{oleg.viro@@gmail.com}
\subjclass{14G30, 14H99; Secondary 14H20, 14N10}

\maketitle

\section*{Introduction}\label{int}

This paper is a translation of the first chapter of my dissertation
\footnote{This is not a Ph D., but a dissertation for the degree of Doctor
of Physico-Mathematical Sciences. In Russia there are two
degrees in mathematics. The lower, degree corresponding approximately
to Ph D., is called Candidate of Physico-Mathematical Sciences.
The high degree dissertation is supposed to be devoted to a subject
distinct from the subject of the Candidate dissertation. My Candidate
dissertation was on interpretation of signature invariants of knots
in terms of intersection form of branched covering spaces of the 4-ball.
It was defended in 1974.} which was defended in 1983.
I do not take here an attempt of updating.

The results of the dissertation were obtained in 1978-80, announced in
\cite{6d,8d,11d}, a short fragment  was published in detail in
\cite{13d} and a considerable part was published in paper \cite{12d}.
The later publication appeared, however, in almost inaccessible edition
and has not been translated into English.

In \cite{Viro:contr.constr.} I presented almost all constructions of
plane curves contained in the dissertation, but in a simplified version:
without description of the main underlying patchwork construction of
algebraic hypersurfaces. Now I regard the latter as the most
important result of the dissertation with potential range of
application much wider than topology of real algebraic varieties. It
was the subject of the first chapter of the dissertation, and it is
this chapter that is presented in this paper.

In the dissertation the patchwork construction
was applied only in the case of plane curves.
It is developed in considerably higher generality.
This is motivated not only by a hope on future
applications, but mainly internal logic of the subject.
In particular, the proof of Main Patchwork Theorem in the
two-dimensional situation is based on results related to the
three-dimensional situation and analogous to the two-dimensional
results which are involved into formulation of the two-dimensional
Patchwork Theorem. Thus, it is natural to formulate and
prove these results once for all dimensions, but then it is not natural
to confine Patchwork Theorem itself to the two-dimensional
case.  The exposition becomes heavier because of high degree of
generality.  Therefore the main text is prefaced with a section
with visualizable presentation of results.
The other sections formally are
not based on the first one and contain the most general formulations
and complete proofs.

In the last section another, more elementary, approach is expounded.
It gives more detailed information about the constructed manifolds,
having not only topological but also metric character. There, in
particular, Main Patchwork Theorem is proved once again.

I am grateful to Julia Viro who translated this text.

\tableofcontents
\pagenumbering{arabic}
\setcounter{section}0

\section{Patchworking plane real algebraic curves }\label{s1}

This Section is introductory. I explain the character of results
staying in the framework of plane curves. A real exposition begins
in Section \ref{s2}. It does not depend on Section \ref{s1}. To a reader
who is motivated enough and does not like informal texts without proofs,
I would recommend to skip this Section.

\subsection{The case of smallest patches}\label{sTC} We start with the
special case of the patchworking. In this case the patches
are so simple that they do not demand a special care. It purifies the
construction and makes it a straight bridge between combinatorial
geometry and real algebraic geometry.

\begin{prop}[Initial Data]\label{TC.A}
 Let $m$ be a positive integer number {\em[}it
is the degree of the curve under construction{\em]}. Let $\GD$  be the
triangle in ${\bf R}^{2}$ with vertices $(0,0)$, $(m,0)$, $(0,m)$
{\em[}it is a would-be Newton diagram of the equation{\em]}. Let $\mathcal
T$ be a triangulation of $\GD$ whose vertices have integer coordinates.
Let the vertices of $\mathcal T$ be equipped with signs;  the sign
(plus or minus) at the vertex with coordinates $(i,j)$ is denoted by
$\Gs_{i,j}$.\end{prop}

See Figure \ref{fTC1}.
\begin{figure}[h]
\centerline{\epsffile{pw-tc1.eps}}
\caption{}
\label{fTC1}
\end{figure}

For $\Ge,\delta=\pm1$ denote the reflection
$\RR^2\to\RR^2:(x,y)\mapsto(\Ge x,\delta y)$ by $S_{\Ge,\delta}$.  For
a set $A\subset\R^2$, denote $S_{\Ge,\delta}(A)$ by
$A_{\Ge,\delta}$ (see Figure \ref{f5}). Denote a quadrant $\{(x,y)\in
\R^2\,|\,\Ge x>0, \delta y>0\}$ by $Q_{\Ge,\delta}$.
\begin{figure}[h]
\centerline{\epsffile{pw-f5wos.eps}} 
\caption{}
\label{f5}
\end{figure}

The following construction associates with Initial Data
\ref{TC.A} above a piecewise linear curve in the projective plane.

\begin{prop}[Combinatorial patchworking] \label{TC.B}
Take the square $\GD_*$ made of
$\GD$ and its mirror images $\GD_{+-}$, $\GD_{-+}$ and $\GD_{--}$.
Extend the triangulation $\mathbb T$ of $\GD$ to a  triangulation
$\mathbb T_*$ of $\GD_*$
symmetric with respect to the coordinate axes. Extend the distribution
of signs $\Gs_{i,j}$ to a distribution of signs on the vertices of the
extended triangulation which satisfies the following condition:
$\Gs_{i,j}\Gs_{\Ge i,\Gd j}\Ge^i\Gd^j=1$ for any vertex $(i,j)$ of
$\mathbb T$ and $\Ge,\delta=\pm1$. {\em(In other words, passing from a
vertex to its mirror image with respect to an axis we preserve its sign
if the distance from the vertex to the axis is even, and change the
sign if the distance is odd.)\footnote{More sophisticated description:
the new distribution should satisfy the modular property:  $g^*
(\Gs_{i,j}x^i y^j) = \Gs_{g(i,j)}x^i y^j$ for $g = S_{\Ge\Gd}$ (in
other words, the sign at a vertex is the sign of the corresponding
monomial in the quadrant containing the vertex).}}

If a triangle of the triangulation $\mathbb T_*$ has vertices
of different signs, draw the midline separating the vertices of
different signs.  Denote by $L$ the union of these midlines. It
is a collection of polygonal lines contained in $\GD_*$.
Glue by $S_{--}$ the opposite sides of $\GD_*$.  The resulting space
$\bar\GD$ is homeomorphic to the projective plane $\rpp$. Denote
by $\bar L$ the image of $L$ in $\bar\GD$. \end{prop}
\begin{figure}[h] \centerline{\epsffile{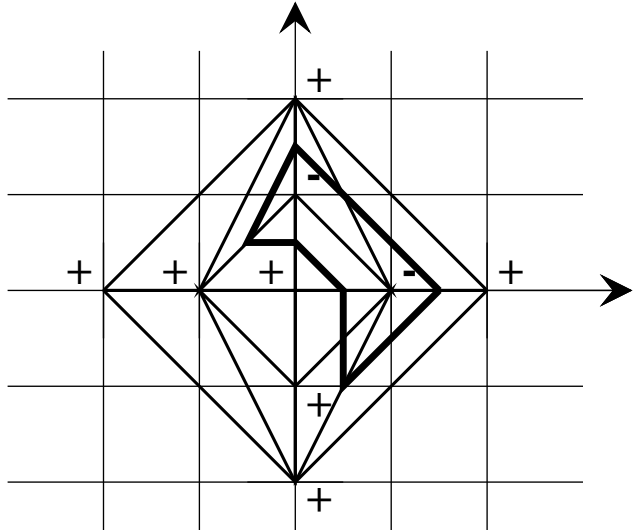}}
\caption{Combinatorial patchworking of the initial data shown in
Figure \ref{fTC1}}
\label{fTC2}
\end{figure}

Let us introduce a supplementary assumption:
the triangulation $\mathbb T$ of $\GD$ is {\it convex}. It means that
there exists a convex piecewise linear function $\nu:\GD \to\R$
which is linear on each triangle of $\mathbb T$ and
not linear on the union of any two triangles of $\mathbb T$. A function
$\nu$ with this property is said to {\it convexify\/} $\mathbb T$.

In fact, to stay in the frameworks of algebraic geometry we need
to accept an additional assumption: a
function $\nu$ convexifying $\mathbb T$ should take integer value on
each vertex of $\mathbb T$. Such a function is said to {\it convexify $\mathbb T$
over\/} $\Z$. However this additional restriction is easy to satisfy.  A
function $\nu:\GD \to\R$ convexifying $\mathbb T$ is characterized by its
values on vertices of $\mathbb T$.  It is easy to see that this provides a
natural identification of the set of functions convexifying $\mathbb T$
with an open convex cone of $\R^N$ where $N$ is the number of vertices
of $\mathbb T$.  Therefore if this set is not empty, then it contains a
point with rational coordinates, and hence a point with integer
coordinates.

\begin{prop}[Polynomial patchworking]\label{TC.C}
 Given Initial Data $m$, $\GD$,
$\mathbb T$ and $\Gs_{i,j}$ as above and a function $\nu$ convexifying
$\mathbb T$ over $\Z$. Take the polynomial
$$b(x,y,t)=\sum_{{\small\begin{aligned}&(i,j)
\text{\enspace runs over}\\ &\text{\quad vertices of $\mathbb
T$}\end{aligned}}}\Gs_{i,j}x^iy^jt^{\nu(i,j)}.  $$ and consider it as a
one-parameter family of polynomials: set $b_t(x,y)=b(x,y,t)$. Denote by
$B_t$ the corresponding homogeneous polynomials:
$$B_t(x_0,x_1,x_2)=x_0^mb_t(x_1/x_0,x_2/x_0).$$\end{prop}

\begin{prop}[Patchwork Theorem]\label{TC.D} Let
$m$, $\GD$, $\mathbb T$ and $\Gs_{i,j}$ be an initial data as above
and $\nu$ a function
convexifying $\mathbb T$ over $\Z$. Denote by $b_t$ and $B_t$ the
non-homogeneous and homogeneous polynomials obtained by the polynomial
patchworking of these initial data and by $L$ and $\bar L$ the piecewise
linear curves in the square $\GD_*$ and its quotient space $\bar\GD$
respectively obtained from the same initial data by the combinatorial
patchworking.

Then there exists $t_0>0$ such that
for any $t\in(0,t_0]$ the equation $b_t(x,y)=0$ defines in the plane
$\R^2$ a curve $c_t$ such that the pair
$(\R^2, c_t)$ is homeomorphic to the pair $(\GD_*,L)$
and the
equation $B_t(x_0,x_1,x_2)=0$
defines in the real projective plane a curve $C_t$ such
that the pair $(\rpp, C_t)$ is homeomorphic to the pair $(\bar\GD,\bar L)$.
\end{prop}

\begin{figure}[t] \centerline{\epsffile{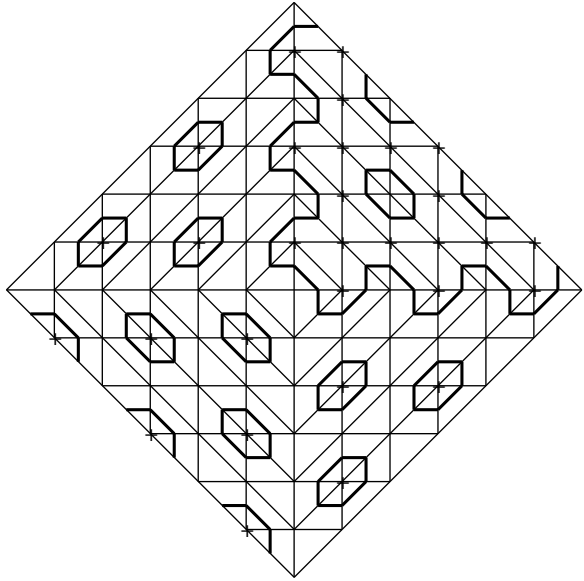}}
\caption{Harnack's curve of degree 6.} \label{fharn} \end{figure}

\begin{figure}[t]
\centerline{\epsffile{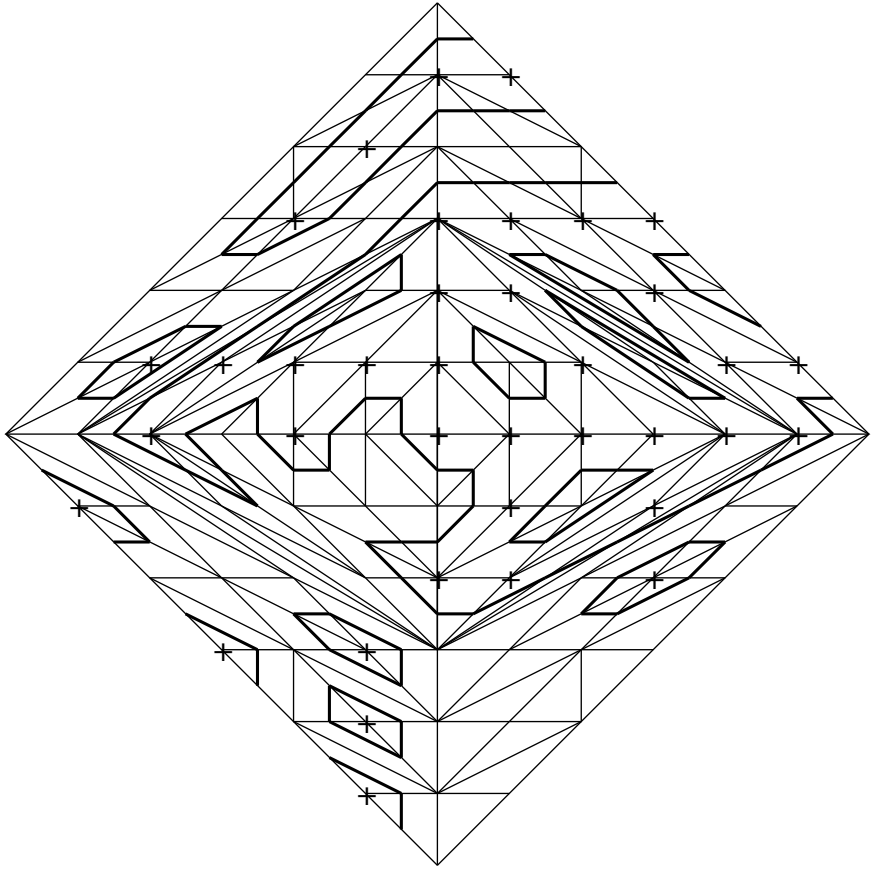}} \caption{Gudkov's
curve of degree 6.} \label{fTC3} \end{figure}

\begin{exmpl}\label{TC.E}
Construction of a curve of degree 2 is shown in Figure
\ref{fTC2}. The broken line corresponds to an ellipse.
More complicated examples with a curves of degree 6 are
shown in Figures \ref{fharn}, \ref{fTC3}.  \end{exmpl}


For more general version of the patchworking we have to prepare
patches. Shortly speaking, the role of patches was played above by
lines. The generalization below is a transition from lines to curves.
Therefore we proceed with a preliminary study of curves.

\subsection{Logarithmic asymptotes of a curve}\label{s1.1} As is
known since Newton's works (see \cite{16s}), behavior of a curve
$\{(x,y)\in\R^2\,|\,a(x,y)=0\}$ near the coordinate axes and at
infinity depends, as a rule, on the coefficients of $a$ corresponding
to the boundary points of its Newton polygon $\Delta(a)$. We need
more specific formulations, but prior to that we have to introduce
several notations and discuss some notions.

For a set $\Gamma\subset\R^2$ and a polynomial
$a(x,y)=\sum_{\omega\in\Bbb Z^2}a_\omega x^{\omega_1}y^{\omega_2}$,
denote the polynomial
$\sum_{\omega\in\Gamma\cap\Bbb Z^2}a_\omega x^{\omega_1}y^{\omega_2}$
by $a^\Gamma$. It is called the {\it $\GG$-truncation\/} of $a$.

For a set $U\subset\R^2$ and a real polynomial $a$ in two variables,
denote  the curve $\{(x,y)\in U\,|\,a(x,y)=0\}$ by  $V_U(a)$.

The complement of the coordinate axes in $\R^2$, i.e. a set
$\{(x,y)\in \R^2\,|\, xy\not=0\}$, is denoted\footnote{This
notation is motivated in Section \ref{s2.3} below.} by $\R\R^2$.

Denote by $l$ the map $\R\R^2\to\R^2$ defined by formula
$l(x,$ $y)=(\ln|x|,\ln|y|)$.
It is clear that the restriction of $l$ to each quadrant
is a diffeomorphism.

A  polynomial in two variables is said to be   {\it
quasi-homogeneous\/} if its Newton polygon is a segment. The
simplest real quasi-homogeneous polynomials are binomials of the form
$\alpha x^p+\beta y^q$ where $p$ and $q$ are relatively prime. A curve
$V_{\Bbb{RR}^2}(a)$, where $a$ is a binomial, is called  {\it
quasiline\/}. The map $l$ transforms quasilines to lines.
In that way any  line with rational slope can be
obtained. The image $l(V_{\Bbb{RR}^2}(a))$ of
the quasiline $V_{\Bbb{RR}^2}(a)$ is orthogonal to the segment
$\Delta(a)$.

It is clear that any real quasi-homogeneous polynomial in 2 variables
is decomposable into a product of binomials of the type described
above and trinomials without zeros in $\Bbb{RR}^2$. Thus if $a$ is a
real quasi-homogeneous polynomial then the curve $V_{\Bbb{RR}^2(a)}$ is
decomposable into a union of several quasilines which are transformed
by $l$ to lines orthogonal to $\Delta(a)$.

A real polynomial $a$ in two variables  is said to be {\it
peripherally nondegenerate} if for any side $\Gamma$ of its Newton
polygon the curve $V_{\Bbb{RR}^2}(a^{\Gamma})$ is nonsingular (it is a
union of quasilines transformed by $l$ to parallel lines, so the
condition that it is nonsingular means absence of multiple
components).  Being peripherally nondegenerate is typical in the
sense that among polynomials with the same Newton polygons the
peripherally nondegenerate ones form  nonempty set
open in the Zarisky topology.

For a side $\Gamma$ of a polygon $\Delta$, denote by
$DC_{\Delta}^-(\Gamma)$ a ray consisting of vectors orthogonal to
$\Gamma$ and directed outside $\Delta$ with respect to $\Gamma$ (see
Figure \ref{f1} and Section \ref{s2.2}).

\begin{figure}[t]
\centerline{\epsffile{pw-f1wos.eps}} 
\caption{}
\label{f1}
\end{figure}

The assertion in the beginning of this Section about behavior of a
curve nearby the coordinate axes and at infinity can be made now more
precise in the following way.

\begin{prop}\label{1.1.A}
Let $\Delta\in\Bbb{RR}^2$ be a convex polygon with
nonempty interior and sides $\Gamma_1,$ \dots, $\Gamma_n$. Let $a$ be
a peripherally nondegenerate real  polynomial in 2 variables with
$\Delta(a)=\Delta$. Then for any quadrant $U\in\Bbb{RR}^2$ each
line contained in $l(V_U(a^{\Gamma_i})$ with $i=1$,\dots,$n$
is an asymptote of $l(V_U(a))$, and $l(V_U(a))$ goes to infinity only
along these asymptotes in the directions defined by rays
$DC^-_\Delta(\Gamma_i)$.
\end{prop}

Theorem generalizing this proposition is formulated in Section
\ref{s6.3} and proved in Section \ref{s6.4}. Here we restrict ourselves
to the following elementary example illustrating \ref{1.1.A}.

\begin{exmpl}\label{1.1.B} Consider the polynomial
$a(x,y)=8x^3-x^2+4y^2$.  Its Newton polygon is shown in Figure
\ref{f1}. In Figure \ref{f2} the curve $V_{\R^2}(a)$ is shown.  In
Figure \ref{f3} the rays $DC^-_\Delta(\Gamma _i)$ and the images of
$V_U(a)$ and $V_U(a^{\Gamma _i})$ under diffeomorphisms $l|_U:U\to\Bbb
R^2$ are shown, where $U$ runs over the set of components of
$\Bbb{RR}^2$ (i.e.  quadrants). In Figure \ref{f4} the images of
$DC^-_\Delta(\Gamma_i)$ under $l$ and the curves $V_{\R^2}(a)$ and
$V_\R^2(a^{\Gamma_i})$ are shown.
\end{exmpl}

\begin{figure}[t]
\centerline{\epsffile{pw-f2.eps}}
\caption{}
\label{f2}
\end{figure}

\begin{figure}[t]
\centerline{\epsffile{pw-f3wos.eps}} 
\caption{}
\label{f3}
\end{figure}

\begin{figure}[b]
\centerline{\epsffile{pw-f4.eps}}
\caption{}
\label{f4}
\end{figure}

\subsection{Charts of polynomials}\label{s1.2} The notion of a chart
of a polynomial is fundamental for what follows. It is introduced
naturally via the theory of toric varieties (see Section \ref{s3}).
Another definition, which is less natural and applicable to a narrower
class of polynomials, but more elementary, can be extracted from the
results generalizing Theorem \ref{1.1.A} (see Section \ref{s6}). In
this Section, first, I try to give a rough idea about the definition
related with toric varieties, and then I give  the definitions related
with Theorem \ref{1.1.A} with all details.

To any convex  closed polygon $\Delta\subset\R^2$ with vertices whose
coordinates are integers, a real algebraic surface $\R\Delta$ is
associated.  This surface is a completion of $\Bbb{RR}^2$
($=(\R\smallsetminus0)^2$). The complement
$\R\Delta \smallsetminus \R\R^2$ consists
of lines corresponding to  sides of $\Delta$. From
the topological viewpoint $\R\Delta$ can be obtained from four
copies of $\Delta$ by pairwise gluing of their sides. For a real
polynomial $a$ in two variables we denote  the closure of
$V_{\R\R^2}(a)$ in $\R\Delta$ by $V_{\R\Delta}(a)$.
Let $a$ be a real polynomial in two variables which is not
quasi-homogeneous. (The
latter assumption is not necessary, it is made for the sake of
simplicity.) Cut the surface $\R\Delta(a)$ along lines of
$\R\Delta(a)\smallsetminus \R\R^2$ (i.e.
replace each of these
lines by two lines). The result is four copies of $\Delta(a)$ and a
curve lying in them obtained from $V_{\R\Delta(a)}(a)$. The pair
consisting of these four polygons and this curve is a chart of $a$.

Recall that for $\Ge,\delta=\pm1$ we denote the reflection
$\RR^2\to\RR^2:(x,y)\mapsto(\Ge x,\delta y)$ by $S_{\Ge,\delta}$.  For
a set $A\subset\R^2$ we denote $S_{\Ge,\delta}(A)$ by
$A_{\Ge,\delta}$ (see Figure \ref{f5}). Denote a quadrant $\{(x,y)\in
\R^2\,|\,\Ge x>0, \delta y>0\}$ by $Q_{\Ge,\delta}$.

Now define the charts for two classes of real polynomials separately.

First, consider the case of quasi-homogeneous polynomials. Let
$a$ be a quasi-homogeneous polynomial defining a nonsingular curve
$V_{\R\R^2}(a)$. Let $(w_1,w_2)$ be a vector orthogonal  to
$\Delta= \Delta(a)$  with integer
relatively prime coordinates. It is
clear that in this case $V_{\R^2}(a)$ is invariant under
$S_{(-1)^{w_1},(-1)^{w_2}}$. A pair $(\Delta_*,$ $\upsilon)$
consisting of $\Delta_*$ and a finite set $\upsilon\subset\Delta_*$
is called  {\it a chart\/} of $a$, if the number of points of
$\upsilon\cap\Delta_{\Ge,\delta}$ is equal to the number of components
of $V_{Q_{\Ge,\delta}}(a)$ and  $\upsilon$ is invariant under
$S_{(-1)^{w_1},(-1)^{w_2}}$ (remind that $V_{\R^2}(a)$ is
invariant under the same reflection).

\begin{exmpl}\label{1.2.A}
In Figure \ref{f6} it is shown a curve
$V_{\R^2}(a)$ with $a(x,y) = 2x^6y - x^4y^2 - 2x^2y^3 + y^4 =
(x^2 - y)(x^2 + y)(2x^2 - y)y$, and  a
chart of $a$.
\begin{figure}[t]
\centerline{\epsffile{pw-f6.eps}}
\caption{}
\label{f6}
\end{figure}
Now consider the case of peripherally nondegenerate polynomials
with Newton polygons having nonempty interiors. Let $\Delta$,
$\Gamma_1,\dots,\Gamma_n$ and $a$ be as in \ref{1.1.A}. Then, as it
follows from \ref{1.1.A}, there exist a disk $D\subset\R^2$ with
center at the origin and neighborhoods $D_1,\dots,D_n$ of rays
$DC^-_\Delta(\Gamma_1),\dots,DC^-_\Delta(\Gamma_n)$ such that the curve
$V_{\R\R^2}(a)$ lies in $l^{-1}(D\cup D_1\cup\dots\cup D_n)$
and for $i = 1,\dots,n$ the curve $V_{l^{-1}(D_i\smallsetminus D)}(a)$
is approximated by $V_{l^{-1}(D_i\smallsetminus D)}(a^{\Gamma_i})$
and can be contracted (in itself) to $V_{l^{-1}(D_i\cap\p D)}(a)$.

A pair $(\Delta_*,$ $\upsilon)$ consisting of $\Delta_*$ and a curve
$\upsilon\subset\Delta_*$ is called a {\it chart} of $a$ if
\begin{enumerate}
\item  for $i = 1,\dots,n$ the pair
$(\Gamma_{i*},$ $\Gamma_{i*}\cap\upsilon)$ is a chart of $a^{\Gamma_i}$
and
\item for $\Ge$, $\delta = \pm1$ there exists a homeomorphism
$h_{\Ge,\delta}:D\to\Delta$ such that
$\upsilon\cap\Delta_{\Ge, \delta} = S_{\Ge,\delta}\circ
h_{\Ge,\delta}\circ l(V_{l^{-1}(D)\cap Q_{\Ge,\delta}}(a))$
and $h_{\Ge,\delta} (\p D\cap D_i)\subset \Gamma_i$ for $i =
1,\dots,n$.  \end{enumerate}

It follows from \ref{1.1.A} that any peripherally nondegenerate real
polynomial $a$ with $\Int\Delta(a)\ne \empt$ has a chart.
It is easy to see that the chart is unique up to a homeomorphism
$\Delta_*\to\Delta_*$ preserving the polygons $\Delta_{\Ge,\delta}$,
their sides and their vertices.
\end{exmpl}

\begin{exmpl}\label{1.2.B}
 In Figure \ref{f8} it is shown a chart of
$8x^3 - x^2 + 4y^2$ which was considered in
\ref{1.1.B}.
\end{exmpl}
\begin{figure}[t]
\centerline{\epsffile{pw-f8.eps}}
\caption{}
\label{f8}
\end{figure}

\begin{prop}[Generalization of Example \ref{1.2.B}]\label{1.2.BG} Let
$$a(x,y)=a_1x^{i_1}y^{j_1}+a_2x^{i_2}y^{j_2}+a_3x^{i_3}y^{j_3}$$
be a non-quasi-homogeneous real polynomial (i.~e., a real trinomial
whose the Newton polygon has nonempty interior).
For $\Ge,\Gd=\pm1$ set
$$\Gs_{\Ge i_k,\Gd j_k}=sign(a_k\Ge^{i_k}\Gd^{j_k}).$$
Then the pair consisting of $\GD_*$ and the midlines of $\GD_{\Ge,\Gd}$
separating the vertices $(\Ge i_k,\Gd j_k)$ with opposite signs
$\Gs_{\Ge i_k,\Gd j_k}$ is a chart of $a$.
\end{prop}

\begin{proof} Consider the restriction of $a$ to the quadrant
$Q_{\Ge,\Gd}$.  If all signs $\Gs_{\Ge i_k,\Gd j_k}$ are the same, then
$a{Q_{\Ge,\Gd}}$ is a sum of three monomials taking values of the same
sign on $Q_{\Ge,\Gd}$. In this case  $V_{Q_{\Ge,\Gd}}(a)$ is empty.
Otherwise, consider the side $\GG$ of the triangle $\GD$ on whose end
points the signs coincide. Take a vector $(w_1,w_2)$ orthogonal to
$\GG$. Consider the curve defined by parametric equation
$t\mapsto (x_0t^{w_1},y_0t^{w_2})$. It is easy to see that the ratio of
the monomials corresponding to the end points of $\GG$ does
not change along this curve, and hence the sum of them is monotone. The
ratio of each of these two monomials with the third one changes from
$0$ to $-\infty$ monotonically. Therefore the trinomial divided by the
monomial which does not sit on $\GG$
changes from $-\infty$ to $1$ continuously and
monotonically. Therefore it takes the zero value once. Curves $t\mapsto
(x_0t^{w_1},y_0t^{w_2})$ are disjoint and fill $Q_{\Ge,\Gd}$.
Therefore, the curve $V_{Q_{\Ge,\Gd}}(a)$ is isotopic to the preimage under
$S_{\Ge,\Gd}\circ h_{\Ge,\Gd}\circ l$
of the midline of the triangle $\GD_{\Ge,\Gd}$ separating the
vertices with opposite signs.\end{proof}

\begin{prop}\label{1.2.C}
 If $a$ is a peripherally nondegenerate real
polynomial in two variables then the topology of a curve
$V_{\R\R^2}(a)$ (i.e. the topological type of pair
$(\R\R^2,V_{\R\R^2}(a))$) and the topology of its
closure in $\R^2$, $\rpp$ and other toric extensions of
$\R\R^2$ can be recovered from a
chart of $a$.  \end{prop}

The part of this proposition concerning to $V_{\R\R^2}(a)$
follows from \ref{1.1.A}. See below Sections
\ref{s2} and \ref{s3} about toric extensions of $\R\R^2$ and
closures of $V_{\R\R^2}(a)$ in them.  In the next Subsection
algorithms recovering the topology of closures of
$V_{\R\R^2}(a)$ in $\R^2$ and $\rpp$ from a chart of $a$ are
described.

\subsection{Recovering the topology of a curve from a
chart of the polynomial}\label{s1.3}First, I shall describe an auxiliary
algorithm which is a block of two main algorithms of this Section.

\begin{prop}[Algorithm. Adjoining a side with normal vector
$(\alpha,\beta)$]\label{1.3.A}
{\em Initial data:} a chart $(\Delta_*,\,\upsilon)$ of a
polynomial.

If $\Delta$ $(=\Delta_{++})$ has a side $\Gamma$ with
$(\alpha,\beta)\in DC^-_\Delta(\Gamma)$ then the algorithm
does not change $(\Delta_*,\,\upsilon)$. Otherwise:

1. Drawn the lines of support  of $\Delta$ orthogonal to
$(\alpha,\beta)$.

2. Take the point belonging to $\Delta$ on each of the two
lines of support, and join these points with a segment.

3. Cut the polygon $\Delta$ along this segment.

4. Move the pieces  obtained  aside from each other by parallel
translations defined by vectors whose difference is orthogonal to
$(\alpha,\beta)$.

5. Fill the space obtained between the pieces with a
parallelogram whose opposite sides are the edges of the cut.

6. Extend the operations applied above to $\Delta$ to
$\Delta_*$ using symmetries $S_{\Ge, \delta}$.

7. Connect the points of edges of the cut obtained from points of
$\upsilon$ with segments which are parallel to the other
pairs of the sides of the parallelograms inserted, and adjoin these
segments to what is obtained from $\upsilon$. The result
and the polygon obtained from $\Delta_*$ constitute the chart
produced by the algorithm.  \end{prop}

\begin{exmpl}\label{1.3.B}
 In Figure \ref{f9}   the steps of Algorithm \ref{1.3.A} are shown. It
is applied to $(\alpha,\beta) = (-1,0)$ and the chart of
$8x^3 - x^2 + 4y^2$ shown in Figure \ref{f8}.
\end{exmpl}
\begin{figure}[t]
\centerline{\epsffile{pw-f9.eps}}
\caption{}
\label{f9}
\end{figure}

Application of Algorithm \ref{1.3.A} to a chart of a polynomial $a$
(in the case when it does change the chart) gives rise a chart of
polynomial $$(x^\beta y^{-\alpha} + x^
{-\beta}y^\alpha)x^{|\beta|}y^{|\alpha|}a(x,y).$$

If $\Delta$ is a segment (i.e. the initial polynomial
is quasi-homogeneous) and this segment is not orthogonal to the vector
$(\alpha,\beta)$ then Algorithm \ref{1.3.A} gives rise to a chart
consisting of four parallelograms, each of which contains as many
parallel segments as components of the curve are contained in
corresponding quadrant.

\begin{prop}[Algorithm] \label{1.3.C}{\sc
Recovering the topology of an affine
curve from a chart of the polynomial.}
{\em Initial data:} a chart $(\Delta_*,\,\upsilon)$ of a polynomial.

1. Apply Algorithm \ref{1.3.A} with $(\alpha,\beta) = (0,-1)$ to
$(\Delta_*,\,\upsilon)$. Assign the former notation
$(\alpha,\beta)$ to the result obtained.

2. Apply Algorithm \ref{1.3.A} with $(\alpha,\beta) = (0,-1)$ to
$(\Delta_*,\,\upsilon)$. Assign the former notation
$(\alpha,\beta)$ to the result obtained.

3. Glue by $S_{+,-}$ the sides of $\Delta_{+,\delta}$, $\Delta_{-,
\delta}$ which are faced to each other and parallel to $(0,1)$
(unless the sides coincide).

4. Glue by $S_{-,+}$ the sides of $\Delta_{\Ge,+}$,
$\Delta_{\Ge,-}$ which are faced to each other and
parallel to $(1,0)$ (unless the sides  coincide).

5. Contract to a point all sides obtained from the sides of
$\Delta$ whose normals are directed into quadrant $P_{-,-}$.

6. Remove the sides which are not touched on in blocks 3, 4 and 5.
\end{prop}

Algorithm \ref{1.3.C} turns the polygon $\Delta_*$
to a space $\Delta'$ which is homeomorphic to $\R^2$, and
the set $\upsilon$  to a set $\upsilon'\subset\Delta'$ such
that the pair $(\Delta',\,\upsilon')$ is homeomorphic to $(\R^2,
\Cl V_{\R\R^2}(a))$, where $\Cl$ denotes closure and $a$ is a
polynomial whose chart is $(\Delta_*,\,\upsilon)$.
\begin{figure}[t]
\centerline{\epsffile{pw-f10.eps}}
\caption{}
\label{f10}
\end{figure}

\begin{exmpl}\label{1.3.D} In Figure \ref{f10}  the steps of Algorithm
\ref{1.3.C} applying to a chart of polynomial $8x^3y - x^2y + 4y^3$ are
shown.  \end{exmpl}

\begin{prop}[Algorithm] \label{1.3.E}{\sc Recovering the topology of
a projective curve from a chart of the polynomial. }
{\em Initial data:} a chart $(\Delta_*,\,\upsilon)$ of a polynomial.

1. Block 1 of Algorithm \ref{1.3.C}.

2. Block 2 of Algorithm \ref{1.3.C}.

3. Apply Algorithm \ref{1.3.A} with $(\alpha,\beta) = (1,1)$ to
$(\Delta_*,\,\upsilon)$. Assign the former notation
$(\Delta_*,\,\upsilon)$ to the result obtained.

4. Block 3 of Algorithm \ref{1.3.C}.

5. Block 4 of Algorithm \ref{1.3.C}.

6. Glue by $S_{-,-}$ the sides of $\Delta_{++}$ and $\Delta_{--}$ which
are faced  to each other and orthogonal to $(1,1)$.

7. Glue by $S_{-,-}$ the sides of $\Delta_{+-}$ and $\Delta_{-+}$ which
are faced to each other and orthogonal to $(1,-1)$.

8. Block 5 of Algorithm \ref{1.3.C}.

9. Contract to a point all sides obtained from the sides of $\Delta$
with normals directed into the angle $\{(x,y)\in\R^2\,|\,x<0,\,
y+x >0\}$.

10. Contract to a point all sides obtained from the sides of $\Delta$
with normals directed into the angle $\{(x,$ $y)\in\R^2\,|\,y <
0,\,y+x >0\}$.
\end{prop}

Algorithm \ref{1.3.E}  turns polygon $\Delta_*$  to
a space $\Delta'$ which is homeomorphic to projective plane $\rpp$,
and the set $\upsilon$  to a set $\upsilon'$ such that the
pair $(\Delta',\,\upsilon')$ is homeomorphic to $(\rpp,\,V_{\R\R^2}(a))$,
where $a$ is the polynomial whose chart is the initial
pair $(\Delta_*,\,\upsilon)$.

\subsection{Patchworking charts}\label{s1.4} Let $a_1,\dots,a_s$ be
peripherally nondegenerate real polynomials in two variables with $\Int
\Delta(a_i)\cap \Int \Delta(a_j) = \empt$ for $i\ne j$. A
pair $(\Delta_*,\,\upsilon)$ is said to be obtained by {\it patchworking\/}
if $\Delta = \bigcup_{i=1}^s\Delta(a_i)$ and there exist charts
$(\Delta(a_i)_*,\,\upsilon_i)$ of $a_1,\dots,a_s$ such that
$\upsilon=\bigcup_{i=1}^s\upsilon_i$.

\begin{exmpl}\label{1.4.A}
In Figure \ref{f8} and Figure \ref{f11} charts
of polynomials $8x^3-x^2+4y^2$ and $4y^2-x^2+1$ are shown. In Figure
\ref{f12} the result of patchworking  these charts is shown.
\end{exmpl}
\begin{figure}[h]
$$\begin{matrix}
{\epsffile{pw-f11.eps}}&\qquad&{\epsffile{pw-f12.eps}}\\
\phantom{A^{A^A}}& & \\
\refstepcounter{figure}\label{f11}
\text{\sc Figure \ref{f11}}& \qquad &
\refstepcounter{figure}\label{f12}
\text{\sc Figure \ref{f12}}
\end{matrix}$$
\end{figure}

\subsection{Patchworking polynomials}\label{s1.5} Let $a_1,\dots,a_s$
be real polynomials in two variables with $\Int \Delta(a_i)\cap \Int
\Delta(a_j) = \empt$ for $i\ne j$ and $a_i^{\Delta(a_i)\cap
\Delta(a_j)}=a_j^{\Delta(a_i)\cap \Delta(a_j)}$ for any $i,$ $j$.
Suppose the set $\Delta = \bigcup_{i=1}^s\Delta(a_i)$ is convex. Then,
obviously, there exists the unique polynomial $a$ with $\Delta(a)=\GD$
and $a^{\Delta(a_i)}=a_i$ for $i=1,\dots,s$.

Let $\nu:\Delta\to\R$ be a convex function such that:
\begin{enumerate}
\item\label{nu1} restrictions $\nu|_{\Delta(a_i)}$ are  linear;
\item\label{nu2} if the restriction of $\nu$ to an open set is linear
then the set is contained in one of $\Delta(a_i)$;
\item\label{nu3}
$\nu(\Delta\cap\Bbb Z^2)\subset\Bbb Z$.  \end{enumerate}
Then $\nu$ is said to {\it convexify\/} the partition
$\GD(a_1),\dots,\GD(a_s)$ of $\GD$.

If $a(x,y)=\sum_{\omega\in\Bbb Z^2}a_\omega x^{\omega_1}y^{\omega_2}$
then we put
$$b_t(x,y)=\sum_{\omega\in\Bbb Z^2}a_\omega
x^{\omega_1}y^{\omega_2}t^{\nu(\omega_1,\omega_2)}$$
and  say that polynomials $b_t$ are obtained by patchworking
$a_1,\dots,a_s$ by $\nu$.

\begin{exmpl}\label{1.5.A}
 Let $a_1(x,y)=8x^3-x^2+4y^2$,
$a_2(x,y)=4y^2-x^2+1$ and
$$\nu(\omega_1,\omega_2)=\begin{cases} 0, & \text{if
$\quad \omega_1+\omega_2\ge2$} \\ 2-\omega_1-\omega_2, & \text{if
$\quad\omega_1+\omega_2\le2$}.\end{cases}$$
Then $b_t(x,y)=8x^3-x^2+4y^2+t^2$.
\end{exmpl}

\subsection{The Main Patchwork Theorem}\label{s1.6} A real
polynomial $a$ in two variables is said to be
{\it completely nondegenerate} if it is
peripherally nondegenerate (i.e. for any side $\Gamma$ of its Newton
polygon the curve $V_{\R\R^2}(a^\Gamma)$ is nonsingular) and
the curve $V_{\R\R^2}(a)$ is nonsingular.

\begin{prop}\label{1.6.A}
If $a_1, \dots, a_s$ are completely nondegenerate
polynomials satisfying all conditions of Section \ref{s1.5}, and $b_t$
are obtained from them by patchworking by some nonnegative convex
function $\nu$ convexifying $\GD(a_1),\dots,\GD(a_s)$,
then there exists $t_0>0$ such
that for any $t\in(0,t_0]$ the polynomial $b_t$ is completely
nondegenerate and its chart is obtained by patchworking charts of $a_1,
\dots, a_s$. \end{prop}

By \ref{1.2.BG}, Theorem \ref{1.6.A} generalizes Theorem \ref{TC.D}.
Theorem  generalizing Theorem \ref{1.6.A} is proven in Section
\ref{s4.3}.  Here we restrict ourselves to several examples.

\begin{exmpl}\label{1.6.B}
 Polynomial $8x^3-x^2+4y^2+t^2$ with
$t>0$ small enough  has the chart shown in Figure \ref{f12}. See
examples \ref{1.4.A} and \ref{1.5.A}.\end{exmpl}

In the next Section there are a number of considerably  more complicated
examples demonstrating efficiency  of
Theorem \ref{1.6.A} in the topology of real algebraic curves.

\subsection{Construction of M-curves of degree 6}\label{s1.7} One of
central points of the well known 16th Hilbert's problem \cite{51d} is
the problem of isotopy classification of curves of degree 6 consisting
of 11 components (by the Harnack inequality \cite{35s} the number of
components of a curve of degree 6 is at most 11). Hilbert
conjectured that there exist only two isotopy types of such curves.
Namely, the types shown in Figure \ref{f13} (a) and (b). His
conjecture was
disproved by Gudkov \cite{14d} in 1969. Gudkov constructed a curve of
degree 6 with ovals' disposition shown in Figure \ref{f13} (c) and
completed
solution of the problem of isotopy classification of nonsingular curves
of degree 6. In particular, he proved, that any curve of degree 6 with
11 components is isotopic to one of the curves of Figure \ref{f13}.
\begin{figure}[t]
\centerline{\epsffile{pw-f13.eps}}
\caption{}
\label{f13}
\end{figure}

Gudkov proposed twice --- in \cite{10s} and \cite{11s} ---
simplified proofs of realizability of the third isotopy type.
His constructions, however, are essentially more complicated than
the construction described below, which is based on \ref{1.6.A} and besides
gives rise to realization of the other two types, and, after a small
modification, realization of almost all isotopy types of nonsingular
plane projective real algebraic curves of degree 6 (see
\cite{Viro:contr.constr.}).

\begin{figure}[t]
\centerline{\epsffile{pw-f14.eps}}
\caption{}
\label{f14}
\end{figure}

{\bf Construction} In Figure \ref{f14} two curves of degree
6  are shown. Each of them has one singular point at which three nonsingular
branches are second order tangent to each other (i.e. this
singularity belongs to type $J_{10}$ in the Arnold classification
\cite{3s}). The curves of Figure \ref{f14} (a) and (b) are easily
constructed by the Hilbert method \cite{50d}, see in
\cite{Viro:contr.constr.}, Section 4.2.

Choosing in the projective plane various affine coordinate
systems, one obtains various polynomials defining these curves. In
Figures \ref{f15} and \ref{f16} charts of four polynomials
appeared in this way are shown. In Figure \ref{f17}  the results of
patchworking charts of Figures \ref{f15} and \ref{f16}
are shown. All constructions can be done in such a way that Theorem
\ref{1.6.A} (see \cite{Viro:contr.constr.}, Section 4.2) may be applied to the corresponding
polynomials.  It ensures existence of polynomials with charts shown in
Figure \ref{f17}.

\begin{figure}[t]
$$\begin{matrix}
&{\epsffile{pw-f15.eps}}\\
& \refstepcounter{figure}\label{f15}
\text{\sc Figure \ref{f15}}\\
\phantom{A^{A^A}}&  \\
{\epsffile{pw-f16.eps}} &{\epsffile{pw-f17.eps}}\\
\refstepcounter{figure}\label{f16}
\text{\sc Figure \ref{f16}}&
\refstepcounter{figure}\label{f17}
\text{\sc Figure \ref{f17}}
\end{matrix}$$
\end{figure}

\subsection{Behavior of curve $V_{\R\R^2}(b_t)$ as
$t\to0$}\label{s1.8} Let $a_1, \dots, a_s$, $\Delta$ and $\nu$
be as in Section \ref{s1.5}.
Suppose that
polynomials $a_1, \dots, a_s$ are completely nondegenerate
 and $\nu|_{\Delta(a_1)}=0$. According
to Theorem \ref{1.6.A}, the polynomial $b_t$ with
sufficiently small $t>0$  has a chart obtained by patchworking  charts of
$a_1, \dots, a_s$.  Obviously, $b_0=a_1$ since
$\nu|_{\Delta(a_1)}=0$. Thus when $t$ comes  to zero the chart of $a_1$
stays only,
the other charts disappear.

How do the domains containing the pieces of $V_{\R\R^2}(b_t)$
homeomorphic to $V_{\R\R^2}(a_1)$, \dots, $V_{\R\R^2}(a_s)$
behave when $t$ approaches zero? They are moving to the coordinate
axes and  infinity. The closer $t$ to zero, the more place is
occupied by the domain, where $V_{\R\R^2}(b_t)$ is organized as
$V_{\R\R^2}(a_1)$ and is approximated by
it (cf. Section \ref{s6.7}).

It is curious that the family $b_t$ can be changed by a simple geometric
transformation
in such a way that the role of $a_1$
passes to any one of $a_2,\dots,a_s$ or even to $a_k^\Gamma$, where
$\Gamma$ is a side of $\Delta(a_k)$, $k=1, \dots, s$. Indeed, let
$\lambda:\R^2\to\R$ be a linear function,
$\lambda(x,y)=\alpha x+\beta y+\gamma$. Let $\nu'=\nu-\lambda$.
Denote by $b'_t$
the result of patchworking $a_1,\dots,a_s$ by $\nu'$. Denote
by $qh_{(a,b),t}$  the linear transformation
$\R\R^2\to\R\R^2:(x,y)\mapsto(xt^a,yt^b)$.
Then
 $$V_{\R\R^2}(b'_t)=V_{\R\R^2}(b_t\circ qh_
{(-\alpha,-\beta),t})=qh_{(\alpha,\beta),t}V_{\R\R^2}(b_t).$$
Indeed,
$$\begin{aligned}
b'_t(x,y) =&\sum a_\omega
x^{\omega_1}y^{\omega_2}t^{\nu(\omega_1,\omega_2)-\alpha\omega_1-
\beta\omega_2 -\gamma} \\
&= t^{-\gamma}\sum
a_\omega(xt^{-\alpha})^{\omega_1}
(yt^{-\beta})^{\omega_2}t^{\nu(\omega_1,\omega_2)} \\
&= t^{-\gamma}b_t(xt^{-\alpha},yt^{-\beta}) \\
&= t^{-\gamma}b_t\circ qh_{(-\alpha,-\beta),t}(x,y).
\end{aligned}$$

Thus the curves $V_{\R\R^2}(b'_t)$ and $V_{\R\R^2}(b_t)$ are
transformed to each other by a linear transformation.
However  the polynomial $b'_t$ does not tend to $a_1$
as $t\to 0$. For
example, if $\lambda|_{\Delta(a_k)}=\nu|_{\Delta(a_k)}$ then
$\nu'|_{\Delta(a_k)}=0$ and $b'_t\to a_k$. In this case as $t\to 0$,
the domains containing parts of $V_{\R\R^2}(b'_t)$, which are
homeomorphic to $V_{\R\R^2}(a_i)$, with $i\ne k$, run away and
the domain in which $V_{\R\R^2}(b'_t)$ looks like
$V_{\R\R^2}(a_k)$ occupies more and more place. If the set,
where $\nu$ coincides with $\lambda$ (or differs from $\lambda$ by a
constant), is a side $\Gamma$ of $\Delta(a_k)$, then  the
curve $V_{\R\R^2}(b'_t)$ turns to $V_{\R\R^2}(a_k^\Gamma)$
(i.e. collection of quasilines) as $t\to0$
similarly.

The whole picture of evolution of $V_{\R\R^2}(b_t)$ when
$t\to0$ is the following. The fragments which look as $V_{\R\R^2}(a_i)$
with $i=1,\dots,s$ become more and more explicit, but
these fragments are not staying. Each of them is moving away from the
others. The only fragment that is growing without moving
corresponds to the set where $\nu$ is  constant. The other fragments
are moving away from it. From the metric viewpoint some of
them (namely, ones going to the origin and axes) are
contracting, while the others are growing. But in the
logarithmic coordinates, i.e. being transformed by
$l:(x,y)\mapsto(\ln|x|,\ln|y|)$, all the fragments are growing (see
Section \ref{s6.7}). Changing $\nu$ we are applying linear
transformation, which distinguishes one fragment and casts away the
others. The transformation turns our attention to a new piece of the
curve. It is as if we would transfer a magnifying lens from one
fragment of the curve to another.  Naturally, under such a
magnification the other fragments disappear at the moment $t=0$.

\subsection{Patchworking as smoothing of singularities}\label{s1.9} In
the projective plane the passage from curves defined by $b_t$ with
$t>0$ to the curve defined by $b_0$ looks quite differently. Here, the
domains, in which the curve defined by $b_t$ looks like curves defined
by $a_1, \dots, a_s$ are not running away, but pressing more closely to
the points $(1:0:0)$, $(0:1:0)$, $(0:0:1)$ and to the axes joining
them. At $t=0$, they are pressed into the points and
axes.  It means that under the inverse passage (from $t=0$ to $t>0$)
the full or partial smoothing of singularities concentrated at the
points $(1:0:0)$, $(0:1:0)$, $(0:0:1)$  and along coordinate axes
happens.
\begin{figure}[t]
\centerline{\epsffile{pw-f18.eps}}
\caption{}
\label{f18}
\end{figure}
\begin{figure}[t]
\centerline{\epsffile{pw-f19.eps}}
\caption{}
\label{f19}
\end{figure}

\begin{exmpl}\label{1.9.A}
 Let $a_1$, $a_2$ be polynomials of degree 6
with
$a_1^{\Delta(a_1)\cap\Delta(a_2)}=a_2^{\Delta(a_1)\cap\Delta(a_2)}$
and charts shown in Figure \ref{f15} (a) and \ref{f16} (b).
Let $\nu_1$, $\nu_2$ and $\nu_3$ be defined by the following
formulas:
$$\begin{aligned}
\nu_1(\omega_1,\omega_2)&=\begin{cases} 0,&\text{if
$\omega_1+2\omega_2\le6$}\\2(\omega_1+2\omega_2-6),&\text{if
$\omega_1+2\omega_2\ge6$}\end{cases} \\
\nu_2(\omega_1,\omega_2)&=\begin{cases} 6-\omega_1-2\omega_2,\phantom{(6)}&\text{if
$\omega_1+2\omega_2\le6$}\\\omega_1+2\omega_2-6,&\text{if
$\omega_1+2\omega_2\ge6$}\end{cases} \\
\nu_3(\omega_1,\omega_2)&=\begin{cases} 2(6-\omega_1-2\omega_2),&\text{if
$\omega_1+2\omega_2\le6$}\\0,&\text{if
$\omega_1+2\omega_2\ge6$}\end{cases}
\end{aligned}
$$
(note, that $\nu_1$, $\nu_2$ and $\nu_3$ differ from each other by
a linear function). Let $b_t^1$, $b_t^2$ and $b_t^3$ be the results of
patchworking $a_1$, $a_2$ by $\nu_1$, $\nu_2$ and $\nu_3$. By Theorem
\ref{1.6.A} for sufficiently small $t>0$  the polynomials $b_t^1$, $b_t^2$
and $b_t^3$ have the same chart shown in Figure \ref{f17} (ab), but as
$t\to 0$ they go to different polynomials, namely, $a_1$,
$a_1^{\Delta(a_1)\cap\Delta(a_2)}$ and $a_2$.The closure of
$V_{\R\R^2}(b_t^i)$ with $i=1$, 2, 3  in the projective plane (they are
transformed to one another by projective transformations) are
shown in Figure \ref{f18}.  The limiting
projective curves, i.e. the projective closures of
$V_{\R\Bbb R^2}(a_1)$, $V_{\R\R^2} (a_1^{\Delta(a_1)\cap\Delta(a_2)})$,
$V_{\R\R^2}(a_2)$ are shown in Figure \ref{f19}. The curve shown in
 Figure \ref{f19} (b) is
the union of three nonsingular conics which are tangent to each other
in two points. \end{exmpl}

Curves of degree 6 with eleven components of all
three isotopy types can be obtained from this curve by small
perturbations of the type under consideration (cf. Section \ref{s1.7}).
Moreover, as it is proven in \cite{Viro:contr.constr.}, Section 5.1,
nonsingular curves of degree 6 of almost all isotopy types can be
obtained.

\subsection{Evolvings of singularities}\label{s1.10} Let $f$ be a real
polynomial in two variables. (See Section 5, where more general
situation with an analytic function playing the role of $f$ is
considered.) Suppose its Newton polygon $\Delta(f)$ intersects both
coordinate axes (this assumption is equivalent to the assumption that
$V_{\R^2}(f)$ is the closure of $V_{\R\R^2}(f)$). Let the distance from
the origin to $\Delta(f)$ be more than 1 or, equivalently, the curve
$V_{\R^2}(f)$ has  a singularity at the origin. Let this singularity be
isolated.  Denote by $B$ a disk with the center at the origin having
sufficiently small radius such that the pair $(B,V_B(f))$ is
homeomorphic to the cone over its boundary
$(\p B,V_{\p B}(f))$ and the curve $V_{\R^2}(f)$
is transversal to $\p B$ (see \cite {15s}, Theorem 2.10).

Let $f$ be included into a continuous family $f_t$ of
polynomials in two variables: $f=f_0$. Such a family is called
a {\it perturbation\/} of $f$. We shall be interested mainly in
perturbations for which curves $V_{\R^2}(f_t)$ have no singular points
in $B$ when $t$ is in some segment of type $(0,\Ge]$. One says about
such a perturbation that it {\it evolves \/} the singularity of
$V_{\Bbb R^2}(f_t)$ at zero.  If perturbation $f_t$ evolves the
singularity of $V_{\R^2}(f)$ at zero then one can find $t_0>0$ such
that for $t\in(0,t_0]$ the curve $V_{\R^2}(f_t)$ has no singularities
in $B$ and, moreover, is transversal to $\p B$. Obviously, there
exists an isotopy $h_t:B\to B$ with $t_0\in(0,t_0]$ such that
$h_{t_0}=\id$ and $h_t(V_B(f_0))=V_B(f_t)$, so all
pairs $(B, V_B(f_t))$ with $t\in(0,t_0]$ are homeomorphic to
each other.  A family $(B,V_{\R^2}(f_t))$ of pairs with
$t\in(0,t_0]$ is called an {\it evolving \/} of singularity of
$V_{\R^2}(f)$  at zero, or an {\it evolving \/} of germ of $V_{\R^2}(f)$.

Denote by $\Gamma_1, \dots,\Gamma_n$ the sides of Newton polygon
$\Delta(f)$ of the polynomial $f$, faced to the origin. Their union
$\Gamma(f)=\bigcup_{i=1}^n\Gamma_i$ is called the
{\it Newton diagram \/} of $f$.

Suppose the curves $V_{\R\R^2}(f^{\Gamma_i})$ with
$i=1,\dots,n$ are nonsingular. Then, according to Newton \cite{16s},
the curve $V_{\R^2}(f)$ is approximated by the union of
$\Cl V_{\R\R^2}(f^{\Gamma_i})$  with $i=1,\dots,n$ in a sufficiently
small neighborhood of the origin. (This is a local version of Theorem
\ref{1.1.A}; it is, as well as \ref{1.1.A}, a corollary of Theorem
\ref{6.3.A}.) Disk $B$ can be taken so small that $V_{\p B}(f)$ is
close to $\p B\cap V_{\R\R^2}(f^{\Gamma_i})$,
so the number and disposition of these points are defined by charts
$(\Gamma_{i*},\,\upsilon_i)$ of $f^\Gamma_i$.  The union
$(\Gamma(f)_*,\,\upsilon)=
(\bigcup_{i=1}^n\Gamma_{i*},\,\bigcup_{i=1}^n\upsilon_i)$
of these charts is called a {\it chart of germ \/} of $V_{\R^2}(f)$ at
zero. It is a pair consisting of a simple closed polygon $\Gamma(f_*)$,
which is symmetric with respect to the axes and encloses the origin,
and finite set $\upsilon$ lying on it. There is a natural bijection of
this set to $V_{\p B}(f)$, which is extendable to a homeomorphism
$(\Gamma(f)_*,$ $\upsilon)\to (\p B,V_{\p B}(f))$.
Denote this homeomorphism by $g$.

Let $f_t$ be a perturbation of $f$, which evolves the singularity at
the origin. Let $B$, $t_0$ and $h_t$ be as above. It is not difficult
to choose an isotopy $h_t:B\to B$, $t\in(0,t_0]$ such that its
restriction to $\p B$ can be extended to an isotopy
$h'_t:\p B\to\p B$ with $t\in[0,t_0]$ and
$h'_0(V_{\p B}(f_{t_0}))=V_{\p B}(f)$.
A pair $(\Pi,\tau)$ consisting of the polygon $\Pi$
bounded by $\Gamma(f)_*$ and an 1-dimensional subvariety $\tau$
of $\Pi$ is called  a {\it chart of evolving \/}
$(B,V_B(f_t))$, $t\in[0,t_0]$ if there exists a homeomorphism
$\Pi\to B$, mapping $\tau$ to $V_{\p B}(f_{t_0})_*$, whose restriction
$\p\Pi\to\p\Pi$ is the composition
$\Gamma(f)_*@>g>>\p B@>{h'_0}>>\p B$.
It is clear that the boundary
$(\p\Pi,\p\tau)$ of a chart of germ's  evolving is a chart of the
germ.  Also it is clear that if polynomial $f$ is completely
nondegenerate and polygons $\Delta(f_t)$ are obtained from $\Delta(f)$
by adjoining the region restricted by the axes and $\Pi(f)$, then
charts of $f_t$ with $t\in(0,t_0]$ can be obtained by patchworking a
chart of $f$ and chart of evolving $(B,V_B(f_t))$, $t\in[0,t_0]$.

The patchworking construction for polynomials gives  a wide class
of evolvings whose charts can be created by the modification of
Theorem \ref{1.6.A} formulated below.

Let $a_1, \dots, a_s$ be completely nondegenerate
polynomials in two variables with
$\Int\Delta(a_i)\cap\Int\Delta(a_j)=\varnothing$ and
$a_i^{\Delta(a_i)\cap\Delta(a_j)}=a_j^{\Delta(a_i)\cap\Delta(a_j)}$ for
$i\ne j$. Let $\bigcup_{i=1}^s\Delta(a_i)$ be a polygon bounded  by
the axes and  $\Gamma(f)$. Let
$a_i^{\Delta(a_i)\cap\Delta(f)}=f^{\Delta(a_i)\cap\Delta(f)}$ for
$i=a, \dots, s$. Let $\nu:\R^2\to\R$ be a nonnegative convex
function which is equal to zero on $\Delta(f)$ and whose restriction on
$\bigcup_{i=1}^s\Delta(a_i)$ satisfies the conditions \ref{nu1},
\ref{nu2} and \ref{nu3} of Section \ref{s1.5} with respect to  $a_1,
\dots,a_s$.  Then a result $f_t$ of patchworking $f$, $a_1, \dots, a_s$
by $\nu$ is a perturbation of $f$.

Theorem \ref{1.6.A} cannot be applied in this situation because the
polynomial $f$ is not supposed to be completely nondegenerate. This
weakening of assumption implies a weakening of conclusion.

\begin{prop}[Local version of Theorem \ref{1.6.A}]\label{1.10.A}
 Under the conditions above
perturbation $f_t$ of $f$ evolves a singularity of $V_{\R^2}(f)$ at
the origin. A chart of the evolving can be obtained by patchworking
charts of $a_1, \dots, a_s$.
\end{prop}

An evolving of a germ, constructed along  the scheme above, is called
a {\it patchwork evolving}.

If  $\Gamma(f)$ consists of one segment and the curve
$V_{\R\R^2}(f^{\Gamma(f)})$ is nonsingular then the germ of
$V_{\R^2}(f)$ at zero is said  to be {\it
semi-quasi-homogeneous}.  In this case for construction of evolving
of the germ of $V_{\R^2}(f)$ according  the scheme above we need  only
one polynomial; by \ref{1.10.A}, its chart  is
a chart of evolving. In this case geometric structure of $V_B(f_t)$
is especially  simple, too: the curve $V_B(f_t)$ is approximated by
$qh_{w,t}(V_{\R^2}(a_1))$,  where $w$ is a vector
orthogonal to $\Gamma(f)$, that is by the curve $V_{\RR^2}(a_1)$ contracted
by the quasihomothety $qh_{w,t}$. Such evolvings were described in my
paper \cite{11d}. It is clear that any patchwork evolving of
semi-quasi-homogeneous germ can be replaced, without changing
its topological models, by a patchwork evolving, in which only one
polynomial is involved (i.e. $s=1$).


\def\Spec{\operatorname{Spec}}
\def\inj{\operatorname{in}}
\section{Toric varieties and their hypersurfaces}\label{s2}

\subsection{Algebraic tori $K\Bbb R^n$}\label{s2.1} In the rest of this
chapter $K$ denotes the main field, which is either  the real
number field $\R$,  or  the complex number field $\C$.

For $\Go=(\Go_1,\dots,\Go_n)\in\Z^n$ and ordered collection
$x$ of variables $x_1, \dots,x_n$ the product
$x_1^{\Go_1}\dots x_n^{\Go_n}$ is denoted by $x^\Go$. A linear
combination of products of this sort
with coefficients from $K$ is called a {\it Laurent polynomial}
or, briefly, {\it L-polynomial} over $K$. Laurent polynomials over
$K$ in $n$ variables form a ring $K[x_1,x_1^{-1},\dots,x_n,x_n^{-1}]$
naturally isomorphic to the ring of regular functions of the
variety $(K\sminus0)^n$.

Below this variety, side by side with the affine space $K^n$ and the
projective space $KP^n$, is one of the main places of action. It is an
algebraic torus over $K$. Denote it by $K\R^n$.

Denote by $l$ the map $K\R^n\to\R^n$ defined by formula
$l(x_1,\dots,x_n)=$ $(\ln|x_1|,$ $\dots,$ $\ln|x_n|)$.

Put $U_K=\{x\in K\,|\,|x|=1\}$, so $U_{\R}=S^0$ and $U_{\C}=S^1$.
Denote by $ar$ the map $K\R^n\to U_K^n$ ($=U_K\times\dots\times U_K)$
defined by
$ar(x_1,\dots,x_n)=(\dfrac{x_1}{|x_1|},\dots,\dfrac{x_n}{|x_n|})$.

Denote by $la$ the map
$$x\mapsto(l(x),ar(x)):K\R^n\to\R^n\times U_K^n.$$
It is clear that this is a diffeomorphism.

$K\R^n$ is a group with respect to the coordinate-wise multiplication,
and $l$, $ar$, $la$ are group homomorphisms; $la$ is an isomorphism of
$K\R^n$ to the direct product of (additive) group $\R^n$ and
(multiplicative) group $U_K^n$.

Being Abelian group, $K\R^n$ acts on itself by translations. Let us
fix notations for some of the translations involved into
this action.

For $w\in\R^n$ and $t>0$ denote by $qh_{w,t}$ and call a {\it
quasi-homothety} with weights $w=(w_1,\dots,w_n)$ and coefficient $t$ the
transformation $K\R^n\to K\R^n$ defined by formula
$qh_{w,t}(x_1,\dots,x_n)=(t^{w_1}x_1,\dots,t^{w_n}x_n)$, i.e. the
translation by $(t^{w_1},\dots,t^{w_n})$. If $w=(1,\dots,1)$ then it
is the usual homothety with coefficient $t$. It is clear that
$qh_{w,t}=qh_{\Gl^{-1}w,t}$ for $\Gl>0$. Denote by $qh_w$ a
quasi-homothety $qh_{w,e}$, where $e$ is the base of natural
logarithms. It is clear, $qh_{w,t}=qh_{(\ln t)w}$.

For $w=(w_1,\dots,w_n)\in U_K^n$ denote by $S_w$ the
translation $K\R^n\to K\R^n$ defined by formula
$$S_w(x_1,\dots,x_n)=(w_1x_1,\dots,w_nx_n),$$
i.~e. the translation by $w$.

For $w\in\R^n$ denote by $T_w$ the translation
$x\mapsto x+w:\R^n\to\R^n$ by the vector $w$.

\begin{prop}\label{2.1.A} Diffeomorphism $la:K\R^n\to\R^n\times U_K^n$
transforms $qh_{w,t}$ to $T_{(\ln t)w}\times \id_{U_K^n}$, and
$S_w$ to $\id_{\R^n}\times(S_w|_{U_K^n})$, i.e.
$$la\circ qh_{w,t}\circ la^{-1}=T_{(\ln t)w}\times \id_{U_K^n}
\quad\text{and}$$
$$la\circ S_w\circ la^{-1}=\id_{\Bbb
R^n}\times(S_w|_{U_K^n}).$$\qed
\end{prop}

In particular, $la\circ qh_w\circ la^{-1}=T_w\times \id$.

A hypersurface of $K\R^n$ defined by $a(x)=0$, where $a$ is a
Laurent polynomial over $K$ in $n$ variables is denoted by
$V_{K\R^n}(a)$.

If $a(x)=\sum_{\Go\in\Z^n}a_\Go x^\Go$ is a Laurent
polynomial, then by its {\it Newton polyhedron} $\GD(a)$ is the
convex hull of $\{\Go\in\R^n\,|\, a_\Go\ne0\}$.

\begin{prop}\label{2.1.B} Let $a$ be a Laurent polynomial over $K$. If
$\GD(a)$ lies in an affine subspace $\GG$ of $\R^n$ then for any
vector $w\in\R^n$ orthogonal to $\GG$, a
hypersurface $V_{K\R^n}(a)$ is invariant under
$qh_{w,t}$.
\end{prop}

\begin{proof} Since $\GD(a)\subset\GG$ and $\GG\perp w$, then
for $\Go\in\GD(a)$ the scalar product $w\Go$ does not depend on
$\Go$. Hence
$$a(qh^{-1}_{w,t}(x))=\sum_{\Go\in\GD(a)}
a_\Go(t^{-w}x)^\Go=t^{-w\Go}\sum_{\Go\in\GD(a)}a_\Go
x^\Go=t^{-w\Go}a(x),$$
and therefore
$$qh_{w,t}(V_{K\R^n}(a))=V_{K\R^n}(a\circ qh_{w,t}^{-1})=
V_{K\R^n}(t^{-w\Go}a)=V_{K\R^n}(a).$$
\end{proof}

Proposition \ref{2.1.B} is equivalent, as it follows from \ref{2.1.A}, to
the assertion that under hypothesis of \ref{2.1.B} the set
$la(V_{K\R^n}(a))$ contains together with each point
$(x,y)\in\R^n\times U_K^n$ all points $(x',y)\in \R^n\times U_K^n$ with
$x'-x\perp\GG$. In other words, in the case $\GD(a)\subset\GG$ the
intersections of $la(V_{K\R^n}(a))$ with fibers $\R^n\times y$ are
cylinders, whose generators are affine spaces of dimension $n-\dim\GG$
orthogonal to $\GG$.

The following proposition can be proven similarly to \ref{2.1.B}.

\begin{prop}\label{2.1.C} Under the hypothesis of \ref{2.1.B} a
hypersurface $V_{K\R^n}(a)$ is invariant under transformations
$S_{(e^{\pi i w_1},\dots,e^{\pi i w_n})}$, where $w\perp\GG$,
$$w\in\begin{cases}\Z^n, \text{\, if $K=\R$}\\
\R^n, \text{\, if $K=\C$.}\end{cases}$$
\qed \end{prop}

In other words, under the hypothesis of \ref{2.1.B} the
hypersurface $V_{K\R^n}(a)$ contains together with each its point
$(x_1,\dots,x_n)$:
\begin{enumerate}
\item points $((-1)^{w_1}x_1,\dots,(-1)^{w_n}x_n)$ with $w\in\Z^n$,
$w\perp\GG$, if $K=\R$,
\item points $(e^{iw_1}x_1,\dots,e^{iw_n}x_n)$ with $w\in\Bbb
R^n$, $w\perp\GG$, if $K=\C$.
\end{enumerate}

\subsection{Polyhedra and cones}\label{s2.2} Below by a {\it
polyhedron} we mean {\it closed convex\/} polyhedron lying in
$\R^n$, which are not necessarily bounded, but have a finite number of
faces. A polyhedron is said to be {\it integer} if on each of its faces
there are enough points with integer coordinates to define the minimal
affine space containing this face. All polyhedra considered below are
assumed to be integer, unless the contrary is stated.

The set of faces of a polyhedron $\GD$ is denoted by $\mathcal{G}(\GD)$,
the set of its $k$-dimensional faces by $\mathcal{G}_k(\GD)$, the set of all
its proper faces by $\mathcal{G}'(\GD)$.

By a {\it halfspace\/} of vector space $V$ we will mean the
preimage of the closed
halfline $\R_+(=\{x\in\R:x\ge0\})$ under a non-zero linear functional
$V\to\R$ (so the boundary hyperplane of a halfspace passes necessarily
through the origin).  By a {\it cone} it is called an intersection of a
finite collection of halfspaces of $\R^n$. A cone is a polyhedron (not
necessarily integer), hence all notions and notations concerning
polyhedra are applicable to cones.

The minimal face of a cone is the maximal vector
subspace contained in the cone. It is called a {\it ridge} of the cone.

For $v_1,\dots,v_k\in\R^n$ denote by $\langle v_1,\dots,v_k\rangle$
the minimal
cone containing $v_1$, \dots, $v_k$; it is called the {\it cone
generated by\/} $v_1,\dots,v_k$. A cone is said to be {\it simplicial}
if it is generated by a collection of linear independent vectors, and
{\it simple} if it is generated by a collection of integer vectors,
which is a basis of the free Abelian group of integer vectors lying in
the minimal vector space which contains the cone.

Let $\GD\subset\R^n$ be a polyhedron and $\GG$ its face.
Denote by $C_\GD(\GG)$ the cone $\bigcup_{r\in\R_+}r\cdot(\GD-y)$,
where $y$ is a point of
$\GG\sminus\partial\GG$. The cone $C_\GD(\GD)$ is
clearly the vector subspace of $\R^n$ which corresponds to the minimal
affine subspace containing $\GD$. The cone $C_\GG(\GG)$ is the
ridge of $C_\GD(\GG)$. If $\GG$ is a face of $\GD$ with
$\dim\GG=\dim\GD-1$, then $C_\GD(\GG)$ is a halfspace of
$C_\GD(\GD)$ with boundary parallel to $\GG$.

For cone $C\subset\R^n$ we put
$$D^+C=\{x\in\R^n\,|\,\forall a\in C \quad ax\ge0\},$$
$$D^-C=\{x\in\R^n\,|\,\forall a\in C \quad ax\le0\}.$$

These are cones, which are said to be {\it dual\/} to $C$. The cones
$D^+C$ and $D^-C$ are symmetric to each other with respect to $0$. The
cone $D^-C$ permits also the following more geometric description. Each
hyperplane of support of $C$ defines a ray consisting of vectors
orthogonal to this hyperplane and directed to that of two open
halfspaces bounded by it, which does not intersect $C$. The union
of all such rays is $D^-C$.

It is clear that $D^+D^+C=C=D^-D^-C$. If $v_1,\dots,v_n$ is a basis of
$\R^n$, then the cone $D^+\langle v_1,\dots,v_n\rangle$ is generated
by dual basis
$v_1^*,\dots,v_n^*$ (which is defined by conditions
$v_i\cdot v_j*=\GD_{ij}$).

\subsection{Affine toric variety}\label{s2.3} Let $\GD\subset\R^n$
be an (integer) cone. Consider the semigroup $K$-algebra $K[\GD\cap\Z^n]$
of the semigroup $\GD\cap\Z^n$. It consists of Laurent
polynomials of the form $\sum_{\Go\in\GD\cap\Z^n}a_\Go x^\Go$.
According to the well known Gordan Lemma (see, for example,
\cite{13s}, 1.3), the semigroup $\GD\cap\Z^n$ is generated by a
finite number of elements and therefore the algebra $K[\GD\cap\Z^n]$
is generated by a finite number of monomials. If this number is greater than
the dimension of $\GD$, then there are nontrivial relations among the
generators; the number of relations of minimal generated collection is
equal to the difference between the number of generators and
the dimension of $\GD$.

{\it An affine toric variety} $K\GD$  is the affine
scheme $\Spec K[\GD\cap\Z^n]$. Its less invariant, but more elementary
definition looks as follows. Let
$$\{\alpha_1,\dots,\alpha_p\,|\,\sum_{i=1}^p u_{1,i}\alpha_i=
\sum_{i=1}^p v_{i,1}\alpha_i,\dots,\sum_{i=1}^p u_{p-n,i}\alpha_i=
\sum_{i=1}^p v_{p-n,i}\alpha_i\}$$
be a presentation of $\GD\cap\Z^n$ by generators and
relations (here $u_{ij}$ and $v_{ij}$ are nonnegative); then the
variety $K\GD$ is isomorphic to the affine subvariety of $K^p$
defined by the system $$\left\{ \begin{aligned} y_1^{u_{11}}\dots
y_p^{u_{1p}}&=y_1^{v_{11}}\dots y_p^{v_{1p}}\\
\hdots\hdots\hdots\hdots&\hdots\hdots\hdots\hdots\\
y_1^{u_{p-n,1}}\dots y_p^{u_{p-n,p}}&=y_1^{v_{p-n,1}}\dots
y_p^{v_{p-n,p}}.
\end{aligned}
\right.
$$
For example, if $\GD=\R^n$, then
$K\GD=\Spec K[x_1,x_1^{-1},\dots,x_n,x_n^{-1}]$ can be presented
as the subvariety of $K^{2n}$ defined by the system
$$\left\{\aligned &y_1y_{n+1}=1\\&\hdots\hdots\hdots\\&y_ny_{2n}=1\endaligned
\right.$$
Projection $K^{2n}\to K^n$ induces an isomorphism of this subvariety to
$(K\sminus0)^n=K\R^n$. This explains the notation $K\R^n$
introduced above.

If $\GD$ is the positive orthant
$A^n=\{x\in\R^n\,|\, x_1\ge0,\dots,x_n\ge0\}$,
 then $K\GD$ is isomorphic to the affine space $K^n$. The same
takes place for any simple cone. If cone is not simple, then
corresponding toric variety is necessarily singular. For example,
the angle shown in Figure \ref{f20} corresponds to  the cone defined in
$K^3$ by $xy=z^2$.

\begin{figure}[t]
\centerline{\epsffile{pw-f20.eps}}
\caption{}
\label{f20}
\end{figure}
Let a cone $\GD_1$ lie in a cone $\GD_2$. Then the inclusion
$\inj:\GD_1\to\GD_2$ defines an inclusion
$K[\GD_1\cap\Z^n]\hookrightarrow K[\GD_2\cap\Z^n]$
which, in turn, defines a regular map
$$\inj^*:\Spec K[\GD_2\cap\Z^n]\to \Spec
K[\GD_1\cap\Z^n],$$
i.e. a regular map $\inj^*:K\GD_2\to
K\GD_1$. The latter can be described in terms of subvarieties of
affine spaces in the following way. The formulas, defining coordinates
of point $\inj^*(y)$ as functions of coordinates of $y$, are
the multiplicative versions
of formulas, defining generators of semigroup $\GD_1\cap\Z^n$
as linear combinations of generators of the ambient
semigroup $\GD_2\cap\Z^n$.

In particular, for any $\GD$ there is a regular map of
$KC_\GD(\GD)\cong K\R^{\dim\GD}$ to $K\GD$. It is not
difficult to prove that it is an open embedding with dense image, thus
$K\GD$ can be considered as a completion of  $K\R^{\dim\GD}$.

An action of algebraic torus $KC_\GD(\GD)$ in itself by
translations is extended to its action in $K\GD$. This extension can
be obtained, for example, in the following way. Note first, that for
defining an action in $K\GD$ it is sufficient to define an action in
the ring $K[\GD\cap\Z^n]$. Define an action of $K\R^n$ on
monomials
$x^\Go\in K[\Gd\cap\ZZ^n]$ by formula
$(\Ga_1,\dots,\Ga_n)x^\Go=\Ga_1^{\Go_1}\dots\Ga_n^{\Go_n}$
and extend it to the whole
ring $K[\GD\cap\Z^n]$ by linearity. Further, note that if
$V\subset\R^n$ is a vector space, then the map $\inj^*:K\R^n\to KV$ is
a group homomorphism. Elements of kernel of
$\inj^*:K\R^n\to KC_\GD(\GD)$ act identically in $K[\GD\cap\Z^n]$. It
allows to extract from the action of $K\R^n$ in $K\GD$ an action of
$KC^\GD(\GD)$ in $K\GD$, which extends the action of
$KC_\GD(\GD)$ in itself by translations.

With each face $\GG$ of a cone $\GD$ one associates (as with a smaller
cone) a variety $K\GG$ and a map $\inj^*:K\GD\to K\GG$.
On the other hand
there exists a map $\inj_*:K\GG\to K\GD$ for which
$\inj^*\circ\,\inj_*$ is
the identity map $K\GG\to K\GG$. Therefore, $\inj_*$ is an embedding
whose
image is a retract of $K\GD$. From the viewpoint of schemes the map
$\inj_*$ should be defined  by the homomorphism
$K[\GD\cap\Z^n]\to K[\GG\cap\Z^n]$
which maps a Laurent polynomial $\sum_{\Go\in\GD\cap\Z^n}a_\Go x^\Go$
to its $\GG$-truncation $\sum_{\Go\in\GG\cap\Z^n}a_\Go x^\Go$.
In terms of subvarieties of affine space, $K\GG$ is the intersection of
$K\GD$ with the subspace
$y_{i_1}=y_{i_2}=\dots=y_{i_s}=0$, where $y_{i_1},\dots,y_{i_s}$ are
the coordinates corresponding to generators of semigroup
$\GD\cap\Z^n$ which do not lie in $\GG$.

Varieties $\inj_*(K\GG)$ with $\GG\in\mathcal{G}_{\dim\GD-1}(\GD)$ cover
$K\GD\sminus\inj^*(KC_\GD(\GD))$. Images of algebraic tori
$KC_\GG(\GG)$ with $\GG\in\mathcal{G}(\GD)$ under the
composition
$$\begin{CD}KC_\GG(\GG)@> \inj^* >> K\GG@> \inj^* >>
K\GD\end{CD}$$
of embeddings form a partition of $K\GD$, which is a smooth
stratification
of $K\GD$. Closure of the stratum $\inj_*\inj^*(KC_\GG(\GG))$ in
$K\GD$ is $\inj_*(K\GG)$. Below in the cases when it does not lead to
confusion we shall identify $K\GG$ with $\inj_*K\GG$ and
$KC_\GG(\GG)$ with $\inj_*\inj^*KC_\GG(\GG)$ (i.e. we shall
consider  $K\GG$ and $KC_\GG(\GG)$ as lying in $K\GD$).

\subsection{Quasi-projective toric variety}\label{s2.4} Let
$\GD\subset\R^n$ be a polyhedron. If $\GG$ is its face and
$\GS$ is a face of
$\GG$, then  $C_\GG(\GS)$ is a face of
$C_\GD(\GG)$ parallel to $\GG$, and
$C_{C_\GD(\GS)}(C_\GG(\GS))=C_\GD(\GG)$, see Figure \ref{f21}.
\begin{figure}[t]
\centerline{\epsffile{pw-f21ws.eps}}
\caption{}
\label{f21}
\end{figure}
In particular, $C_\GD(\GS)\subset C_\GD(\GG)$ and, hence,
the map $\inj^*:KC_\GD(\GG)\to KC_\GD(\GS)$ is defined. It is
easy to see that this is an open embedding. Let us glue all
$KC_\GD(\GG)$ with $\GG\in\mathcal{G}(\GD)$ together by these
embeddings. The result is denoted by $K\GD$ and called the {\it toric
variety} associated with $\GD$. This definition agrees with the
corresponding
definition from the previous Section: if $\GD$ is a cone and $\GS$ is
its ridge then $C_\GD(\GS)=\GD$ and, since the ridge is the minimal
face, all $KC_\GD(\GG)$ with $\GG\in\mathcal{G}(\GD)$ are embedded in
$KC_\GD(\GS)$ and the gluing gives
$KC_\GD(\GS)=K\GD$.

For any polyhedron $\GD$ the toric variety $K\GD$ is quasi-projective.
If $\GD$ is bounded, it is projective (see \cite{54d} and
\cite{13s}).

A polyhedron $\GD\subset\R^n$ is said to be {\it permissible} if
$\dim\GD=n$, each face of $\GD$ has a vertex and for any vertex
$\GG\in\mathcal{G}_0(\GD)$ the cone $C_\GD(\GG)$ is simple. If
polyhedron
$\GD$ is permissible then variety $K\GD$ is nonsingular and it
can be obtained by gluing affine spaces $KC_\GD(\GG)$
with $\GG\in\mathcal{G}_0(\GD)$. The  gluing allows the following
description. Let us associate with each cone $C_\GD(\GG)$ where
$\GG\in\mathcal{G}_0(\GD)$ an automorphism $f_\GG:K\R^n\to K\R^n$:
if $C_\GD(\GG)=\langle v_1,\dots,v_n\rangle$ and
$v_i=(v_{i1},\dots,v_{in})$ for $i=1,\dots,n$, then we put
$f_\GG(x_1,\dots,x_n)=(x_1^{v_{11}}\dots x_n^{v_{1n}},\dots,
x_1^{v_{n1}}\dots x_n^{v_{nn}})$. The variety $K\GD$ is obtained by
gluing to $K\R^n$  copies of $K^n$
by maps
$\begin{CD}K\R^n@> f_\GG>>K\R^n\hookrightarrow K^n\end{CD}$ for all vertices
$\GG$ of $\GD$. (Cf. Khovansky \cite{27s}.)

The variety $K\GD$ is defined by $\GD$, but does not define it.
Indeed, if $\GD_1$ and $\GD_2$ are polyhedra such that there
exists a bijection $\mathcal{G}(\GD_1)\to\mathcal{G}(\GD_2)$, preserving dimensions
and inclusions and assigning to each face of $\GD_1$  a parallel face of
$\GD_2$, then $K\GD_1=K\GD_2$.

Denote by $P^n$ the simplex of dimension $n$ with vertices
$$(0,0,\dots,0),
(1,0,\dots,0), (0,1,0,\dots,0), \dots, (0,0,\dots,1).$$
It is permissible
polyhedron.  $KP^n$ is the n-dimensional projective space (this agrees
with its usual notation).

Evidently, $K(\GD_1\times\GD_2)=K\GD_1\times K\GD_2$. In
particular, if  $\GD\subset\R^2$ is a square with vertices
$(0, 0)$, $(1, 0)$, $(0, 1)$ and $(1, 1)$, i.e. if
$\GD=P^1\times P^1$, then $K\GD$ is a surface  isomorphic to
nonsingular projective surface of degree 2 (to hyperboloid in the case
of $K=\R^2$).

Polyhedra shown in Figure \ref{f22} define the following surfaces:
$K\GD_1$
is the affine plane with a point blown up; $K\GD_2$ is projective plane
with a point blown up ($\R\GD_2$ is the Klein bottle); $K\GD_3$
is the linear surface over $KP^1$, defined by sheaf $\mathcal{O}+\mathcal{O}(-2)$
($\R\GD_3$ is homeomorphic to torus).
\begin{figure}[t]
\centerline{\epsffile{pw-f22ws.eps}}
\caption{}
\label{f22}
\end{figure}

The variety $KC_\GD(\GD)$ is isomorphic to $K\R^{\dim\GD}$,
open and dense in $K\GD$, so $K\GD$ can be considered as a
completion of $K\R^{\dim\GD}$. Actions of $KC_\GD(\GD)$ in
affine parts $KC_\GD(\GG)$ of $K\GD$ correspond to each other
and define an action in $K\GD$ which is an extension of the action of
$KC_\GD(\GD)$ in itself by translations. Transformations of
$K\GD$ extending $qh_{w,t}$ and $S_w$ are denoted by the same
symbols $qh_{w,t}$ and $S_w$.

The complement $K\GD\sminus KC_\GD(\GD)$ is covered by
$K\GS$ with $\GS\in\mathcal{G}(C_\GD(\GG))$,
$\GG\in\mathcal{G}'(\GD)$ or, equivalently, by varieties
$KC_\GG(\GS)$ with $\GS\in\mathcal{G}(\GG)$,
$\GG\in\mathcal{G}'(\GD)$. They comprise  varieties $K\GG$ with
$\GG\in\mathcal{G}'(\GD)$, which also cover
$K\GD\sminus KC_\GD(\GD)$.
The varieties $K\GG$ are situated with respect to  each other in the same
manner as the corresponding faces in the polyhedron:
$K(\GG_1\cap\GG_2)=K\GG_1\cap K\GG_2$. Algebraic tori
$KC_\GG(\GG)=K\GG\sminus\bigcup_{\GS\in\mathcal{G}'(\GG)}K\GS$ form partition
of $K\GD$, which is a smooth stratification; they are orbits of the
action of $KC_\GD(\GD)$ in $K\GD$.

We shall say that a polyhedron $\GD_2$ is {\it richer} than
a polyhedron $\GD_1$ if for any face $\GG_2\in\mathcal{G}(\GD_2)$
there exists a face $\GG_1\in\mathcal{G}(\GD_1)$ such that
$C_{\GD_2}(\GG_2)\supset C_{\GD_1}(\GG_1)$ (such a face
$\GG_1$ is automatically unique), and for each face
$\GG_1\in\mathcal{G}(\GD_1)$ the cone $C_{\GD_1}(\GG_1)$ can be
presented as the
intersection of several cones $C_{\GD_2}(\GG_2)$ with
$\GG_1\in\mathcal{G}(\GD_2)$. This definition allows a convenient
reformulation in terms of dual cones: a polyhedron $\GD_2$ is richer
than polyhedron $\GD_1$ iff the cones
$D^+C_{\GD_2}(\GG_2)$ with $\GG_2\in\mathcal{G}(\GD_2)$ cover
the  set, which is covered by $D^+C_{\GD_1}(\GG_1)$ with
$\GG_1\in\mathcal{G}(\GD_1)$, and the first covering is a refinement
of the second.

Let a polyhedron $\GD_2$ be richer than $\GD_1$. Then the inclusions
$C_{\GD_1}(\GG_1)\hookrightarrow C_{\GD_2}(\GG_2)$ define
for any $\GG_2\in\mathcal{G}(\GD_2)$ a regular map
$\begin{CD}KC_{\GD_2}(\GG_2)@>\inj^*>>KC_{\GD_1}(\GG_1)\hookrightarrow K\GD_1\end{CD}$.
Obviously, these maps commute with the embeddings, by which  $K\GD_2$
and $K\GD_1$ are glued from affine pieces, thus a regular map
$K\GD_2\to K\GD_1$ appears.

One can show (see, for example, \cite{54d}) that for any polyhedron
$\GD_1$ there exists a richer polyhedron $\GD_2$, defining a
nonsingular toric variety $K\GD_2$. Such a polyhedron is called a
{\it resolution\/} of $\GD_1$ (because it  gives a
resolution of singularities of $K\GD_1$). If $\dim\GD=n$ ($=$
the dimension of the ambient space $\R^n$), then a resolution of $\GD$
can be found among permissible polyhedra.

\subsection{Hypersurfaces of toric varieties}\label{s2.5} Let
$\GD\subset\R^n$ be a polyhedron and $a$ be a Laurent polynomial
over $K$ in $n$ variables. Let $C_{\GD(a)}(\GD(a))\subset C_\GD(\GD)$.
Then there exists a monomial $x^\Go$ such that
$\GD(x^\Go a)\subset C_\GD(\GD)$.  The hypersurface $V_{KC_\GD(\GD)}$
does not depend on the
choice of $x^\Go$ and is denoted simply by $V_{KC_\GD(\GD)}(a)$. Its
closure in $K\GD$ is denoted by $V_{K\GD}(a)$. \footnote{Here it is
meant the closure of $K\GD$ in  the Zarisky topology; in the case of
$K=\Bbb C$ the classic topology gives the same result, but in the case
of $K=\Bbb R$ the usual closure may be a nonalgebraic set.} Thus, to
any Laurent polynomial $a$ over $K$ with
$C_{\GD(a)}(\GD(a))\subset C_\GD(\GD)$,  a hypersurface
$V_{K\GD}(a)$ of $K\GD$ is related.
For Laurent polynomial $a(x)=\sum_{\Go\in\Z^n}a_\Go x^\Go$ and
a set $\GG\subset\R^n$  a Laurent
polynomial $a(x)=\sum_{\Go\in\GG\cap\Z^n}a_\Go x^\Go$
is denoted by $a^\GG$ and called the $\GG$-{\it truncation} of $a$.

\begin{prop}\label{2.5.A} Let $\GD\subset\R^n$ be a polyhedron and $a$ be
a Laurent polynomial over $K$ with $C_{\GD(a)}(\GD(a))\subset C_\GD(\GD)$.
If $\GG_1\in\mathcal{G}'(\GD(a))$, $\GG_2\in\mathcal{G}'(\GD)$ and
$C_{\GD(a)}(\GG_1)\subset C_\GD(\GG_2)$ then
$K\GG_2\cap V_{K\GD}(a)=V_{K\GG_2}(a^{\GG_1}).$
\end{prop}

\begin{proof}Consider
$KC_\GD(\GG_2)$. It is a dense subset of $K\GG_2$.
Since $C_{\GD(a)}(\GG_1)\subset C_\GD(\GG_2)$,
there exists a monomial $x^\Go$
such that $\GD(x^\Go a)$ lies in  $C_\GD(\GG_2)$
and intersects its ridge exactly in the face obtained from
$\GG_1$. Since on $K\GG_2\cap KC_\GD(\GG_2)$ all
monomials, whose exponents do not lie on ridge
$C_{\GG_2}(\GG_2)$ of $C_\GD(\GG_2)$, equal zero,
it follows that the intersection
$\{x\in KC_\GD(\GG_2)\, |\, x^\Go a(x)=0\}\cap K\GG_2$ coincides
with
$\{x\in KC_\GD(\GG_2)\, |\,[x^\Go a]^{C_{\GG_2}(\GG_2)}(x)=0\}\cap K\GG_2$.
Note finally, that the latter coincides with
$V_{K\GG_2}(a_1)$.  \end{proof}

\begin{prop}\label{2.5.B}Let $\GD$ and $a$ be as in \ref{2.5.A} and
$\GG_2$  be  a proper face of the polyhedron $\GD$. If there
is no face $\GG_1\in\mathcal{G}'(\GD(a))$ with
$C_{\GD(a)}(\GG_1)\subset C_\GD(\GG_2)$ then
$K\GG_2\subset V_{K\GD}(a)$.
\end{prop}

The proof is analogous to the proof of the previous statement.\qed

Denote by $SV_{K\R^n}(a)$ the set of singular points of
$V_{K\R^n}(a)$, i.e. a set $V_{K\R^n}(a)\cap\bigcap_{i=1}^n
V_{K\R^n}(\frac{\partial a}{\partial x_i})$.

A Laurent polynomial $a$ is said to be {\it completely nondegenerate}
(over $K$) if, for any face $\GG$ of its Newton polyhedron,
$SV_{K\R^n}(a^\GG)$ is empty and, hence, $V_{K\R^n}
(a^\GG)$ is a nonsingular hypersurface. A Laurent polynomial $a$ is
said to be {\it peripherally nondegenerate} if for any proper face
$\GG$ of its Newton polyhedron $SV_{K\R^n}(a^\GG)=\empt$.

It is not difficult to prove that completely nondegenerate
L-polynomials form Zarisky open subset of the space of L-polynomials
over $K$ with a given Newton polyhedron, and the same holds true also
for peripherally nondegenerate L-polynomials.

\begin{prop}\label{2.5.C} If a Laurent polynomial $a$ over $K$ is
completely nondegenerate and  $\GD\subset\R^n$ is a
resolution of its Newton polyhedron $\GD(a)$ then the variety
$V_{K\GD}(a)$ is nonsingular and transversal to all $K\GG$
with $\GG\in\mathcal{G}'(\GD)$. See, for example, \cite{27s}.  \qed\end{prop}

Theorem \ref{2.5.C} allows various  generalizations related with
possibilities to consider singular $K\GD$ or only some
faces of $\GD(a)$ (instead of all of them).
For example, one can show that if under
the hypothesis of \ref{2.5.A} a truncation $a^\GG$ of
$a$ is completely nondegenerate then under an appropriate
understanding  of transversality (in the sense of stratified space
theory) $V_{K\GD}(a)$ is transversal to $K\GG_2$. Without
going into discussion of transversality in this situation, I
formulate a special case of this proposition, generalizing Theorem
\ref{2.5.C}.

\begin{prop}\label{2.5.D} Let $\GG$ be a face of a polyhedron
$\GD\subset\R^n$ with nonempty $\mathcal{G}_0(\GG)$ and with
simple cones $C_\GD(\GS)$ for all $\GS\in\mathcal{G}_0(\GG)$. Let
$a$ be a Laurent polynomial over $K$ in $n$ variables and $\GG_1$ be
a face of $\GD(a)$ with $C_{\GD(a)}(\GG_1)\subset
C_\GD(\GG)$. If  $a^\GG$ is completely
nondegenerate, then the set of singular points of
$V_{K\GD}(a)$ does not intersect $K\GG$ and $V_{K\GD}(a)$
is transversal to $K\GG$.
\end{prop}

The proof of this proposition is a fragment of the proof of Theorem
\ref{2.5.C}.\qed

\begin{prop}[(Corollary of \ref{2.1.B} and \ref{2.1.C})]\label{2.5.E} Let
$\GD$ and $a$ be as in \ref{2.5.A}. Then for any vector $w\in
C_\GD(\GD)$ orthogonal to $C_{\GD(a)}(\GD(a))$, a hypersurface
$V_{K\GD}(a)$ is invariant under transformations $qh_{w,t}:K\GD\to
K\GD$ and $S_{(e^{\pi iw_1},\dots,e^{\pi iw_n)}}:K\GD\to K\GD$
(the  latter in the case of $K=\R$ is defined only if $w\in\Z^n$).
\qed\end{prop}

\def\inj{\operatorname{in}}

\section{Charts}\label{s3}
\subsection{Space $\R_+\GD$}\label{s3.1}The aim of this
Subsection is to
distinguish in $K\GD$ an important  subspace which looks  like $\GD$.
More precisely, it is defined a stratified real
semialgebraic variety $\R_+\GD$, which is embedded in $K\GD$
and homeomorphic, as a stratified space, to the polyhedron $\GD$
stratified by its faces. Briefly $\R_+\GD$ can be described  as the
set of points with nonnegative real coordinates.

If $\GD$ is a cone then $\R_+\GD$ is defined as a subset of
$K\GD$ consisting of the points in which values of all monomials
$x^\Go$ with $\Go\in\GD\cap\Z^n$ are real and nonnegative.
It is clear that for $\GG\in\mathcal{G}'(\GD)$ the set $\mathbb{R}_+
\GG$ coincides with $\R_+\GD\cap K\GG$ and for cones
$\GD_1\subset\GD_2$ a preimage of $\R_+\GD_1$ under
$\inj^*:K\GD_2\to K\GD_1$ (see Section \ref{s2.3}) is $\R_+\GD_2$.

Now let $\GD$ be an arbitrary polyhedron. Embeddings, by which
$K\GD$ is glued form $KC_\GD(\GG)$ with $\GG\in\mathcal{G}(\GD)$, 
embed the sets $\R_+C_\GD(\GG)$ in one another; a
space obtained by gluing from $\R_+C_\GD(\GG)$ with
$\GG\in\mathcal{G}(\GD)$ is $\R_+\GD$. It is clear that if
$\GG\in\mathcal{G}'$ then $\R_+\GG=\R_+\GD\cap K\GG$.

$\R_+\R^n$ is the open positive orthant $\{x\in\Bbb {RR}^n\,|\,
x_1>0,\dots,x_n>0\}$. It can be identified with the
subgroup of quasi-homotheties
of $K\R^n$: one assigns $qh_{l(x)}$ to a point $x\in\Bbb{R_+R}^n$.

If $A^n=\{x\in\R^n | x_1\ge0,\dots,x_n\ge0\}$ then $KA^n=K^n$ (cf.
Section \ref{s2.3}) and $\R_+A^n=A^n$.

If $P^n$ is the $n$-simplex with vertexes $(0,0,\dots,0)$,
$(1,0,\dots,0)$, $(0,1,0, \dots, 0)$, \dots, $(0,0,\dots,1)$,
then $KP^n$ is the $n$-simplex consisting of points of
projective space with nonnegative real homogeneous coordinates.

The set $\R_+\GD$ is invariant under quasi-homotheties. Orbits of
action in $\R_+\GD$ of the group of quasi-homotheties of $\R_+\R^n$
are sets $\R_+C_\GG(\GG)$ with $\GG\in\mathcal{G}(\GD)$. Orbit
$\R_+C_\GG(\GG)$ is homeomorphic to $\R^{\dim\GG}$ or,
equivalently,   to the interior of $\GG$. Closures $\R_+\GG$ of
$\R_+C_\GG(\GG)$ intersect one another in the same manner as the
corresponding faces:  $\R_+\GG_1\cap\R_+\GG_2=\R_+(\GG_1\cap\GG_2)$.
From this and from the fact that $\R_+\GG$ is locally conic
(see \cite{57d}) it
follows that $\R_+\GD$ is homeomorphic, as a stratified space, to
$\GD$. However, there is an explicitly constructed homeomorphism.
It is provided by the Atiyah moment map \cite{72d} and in the case
of bounded $\GD$ can be described in the following way.

Choose a collection of points $\Go_1,\dots,\Go_k$ with integer
coordinates, whose convex hull is $\GD$. Then for $\GG\in\mathcal{G}(\GD)$ 
and $\Go_0\in\GG\smallsetminus\partial\GG$ cone
$C_\GD(\GG)$ is $\langle\Go_1-\Go_0,\dots,\Go_k-\Go_0
\rangle$.
For $y\in KC_\GD(\GG)$ denote by $y^\Go$ a value of monomial
$x^\Go$ where $\Go\in C_\GD(\GG)\cap\Z^n$ at this
point. Put
$$M(y)=\dfrac{\sum_{i=1}^{k}|y^{\Go_i-\Go_0}|\Go_i}
{\sum_{i=1}^{k}{|y^{\Go_i-\Go_0}|}}\in\R^n.$$
Obviously $M(y)$ lies in $\GD$, does non depend on the choice of
$\Go_0$ and for $y\in KC_\GD(\GG_1)\cap KC_\GD(\GG_2)$
does not depend on what face, $\GG_1$ or $\GG_2$, is used for the
definition of $M(y)$. Thus a map $M:K\GD\to\GD$ is well defined. It
is not difficult to show that $M|_{\R_+\GD}:\R_+\GD\to\GD$
is a stratified homeomorphism.

\subsection{Charts of $K\GD$}\label{s3.2}The space $K\R^n$ can be
presented as
$\Bbb{R_+R}^n\times U_K^n$. In this Section an analogous
representation of $K\GD$ is described.

$\R_+\GD$ is a fundamental domain for the natural action of $U_K^n$
in
$K\GD$, i.e. its intersection with each orbit of the action consists
of one point.

For a point $x\in\R_+C_\GG\GG$ where $\GG\in\mathcal{G}(\GD)$, 
the stationary subgroup of action of $U_K^n$ consists of
transformations $S_{(e^{\pi iw_1},\dots,e^{\pi iw_n})}$, where vector
$(w_1,\dots,w_n)$ is orthogonal to $C_\GG(\GG)$. In particular,
if $\dim\GG=n$ then the stationary subgroup is trivial. If
$\dim\GG=n-r$ then it is isomorphic to $U_K^r$. Denote by $U_\GG$
a subgroup of $U_K^n$ consisting of elements $(e^{\pi
iw_1},\dots,e^{\pi iw_n})$ with $(w_1,\dots,w_n)\perp
C_\GG(\GG)$.

Define a map $\rho:\R_+\GD\times U_K^n\to K\GD$ by formula
$(x,y)\mapsto S_y(x)$. It is surjection and we know the partition
of $\R_+\GD\times U_K^n$ into  preimages of points. Since $\rho$ is
proper and $K\GD$ is locally compact and Hausdorff, it follows
that $K\GD$ is homeomorphic to the quotientspace of
$\R_+\GD\times U_K^n$ with respect to the partition into sets
$x\times yU_\GG$ with
$x\in\R_+C_\GG(\GG)$, $ y\in U_K^n$.

Consider as an example the case of $K=\R$ and $n=2$. Let
a polyhedron $\GD$ lies in the open positive quadrant. We place
$\GD\times U_\R^2$ in $\R^2$  identifying $(x,y)\in
\GD\times U_\R^2$ with $S_y(x)\in\R^2$. $\Bbb
R_+\GD\times U_\R^2$ is homeomorphic $\GD\times U_\R^2$,
so the surface $\R\GD$ can be obtained by an appropriate gluing
(namely, by transformations taken from $U_\GG$)  sides of four polygons
consisting $\GD\times U_\R^2$. Figure \ref{f23}  shows what
gluings ought to be done in three special cases.
\begin{figure}[t]
\centerline{\epsffile{pw-f23.eps}}
\caption{}
\label{f23}
\end{figure}

\subsection{Charts of L-polynomials}\label{s3.3} Let $a$ be a Laurent
polynomial over $K$ in $n$ variables and $\GD$ be its Newton
polyhedron.  Let $h$ be a homeomorphism $\GD\to\R_+\GD$,
mapping each face to the corresponding subspace, and such that for any
$\GG\in\mathcal{G}(\GD)$, $x\in\GG$, $y\in U^n_K$, $z\in U_\GG$ \quad
$$h(x,y,z)=(pr_{\R_+\GG}h(x,y),\,zpr_{U_K^n}h(x,y)).$$
For $h$ one can
take, for example, $(M|_{\R_+\GD})^-$.

A pair consisting of $\GD\times U_K^n$ and its subset
$\upsilon$ which is the preimage of $V_{K\GD}(a)$ under
$$\begin{CD}\GD\times U_K^n@>h\times \id>>\R_+\GD\times U_K^n@>\rho>>K\GD
\end{CD}$$
is called a (nonreduced) $K$-{\it chart} of L-polynomial $a$.

It is clear that the set $\upsilon$ is invariant under transformations
$\id\times S$ with $S\in U_\GD$ and its intersection with $\GG\times
U_K^n$, where $\GG\in\mathcal{G}'(\GD)$ is invariant under transformations
$\id\times S$ with $S\in U_\GG$.

As it follows from \ref{2.5.A}, if $\GG$ is a face of $\GD$, and $(\GD\times
U_K^n, \,\upsilon)$ is a nonreduced $K$-chart of L-polynomial $a$, then
$(\GG\times U_K^n, \,\upsilon\cap(\GG\times U_K^n))$ is a nonreduced
$K$-chart of L-polynomial $a^\GG$.

A nonreduced $K$-chart of Laurent polynomial $a$ is unique up to
homeomorphism $\GD\times U_K^n\to\GD\times U_K^n$, satisfying the
following two conditions:
\begin{enumerate}
\item it map $\GG\times y$  with
$\GG\in\mathcal{G}(\GD)$ and $y\in U_K^n$ to itself and
\item its
restriction to $\GG\times U_K^n$ with $g\in\mathcal{G}(\GD)$
commutes  with transformations $\id\times S:\GG\times
U_K^n\to\GG\times U_K^n$ where $S\in U_\GG$.
\end{enumerate}

In the case when $a$ is a usual polynomial, it is convenient to place its
$K$-chart into $K^n$. For this, consider a map $A^n\times U_K^n\to
K^n:(x,y)\mapsto S_y(x)$. Denote by $\GD_K(a)$ the image of
$\GD(a)\times U_K^n$ under this map. Call by a (reduced) $K$-{\it chart}
of $a$ the  image of a nonreduced $K$-chart of $a$ under this map.
The charts of peripherally nondegenerate real polynomial in
two variables introduced in Section \ref{s1.2} are $\R$-charts in the
sense of this definition.

\begin{prop}\label{3.3.A}Let $a$ be a Laurent polynomial over $K$ in $n$
variables, $\GG$  a face of its Newton polyhedron,
$\rho:\R_+\GD(a)\times U_K^n\to K\GD(a)$  a natural projection.
If the truncation $a^\GG$ is completely nondegenerate then the
set of singular points of hypersurface
$\rho^{-1}V_{K\GD(a)}(a)$ of \,$\R_+\GD(a)\times U_K^n$
does not intersect $\R_+\GG\times U_K^n$, and $\rho^{-1}
V_{K\GD(a)}(a)$ is transversal to $\R_+\GG\times U_K^n$.
\end{prop}
\begin{proof}Let $\GD$ be a resolution of polyhedron $\GD(a)$. Then a
commutative diagram
$$\begin{CD}
(\R_+\GD\times U_K^n,\, \rho^{\prime-1}(V_{K\GD}(a)))
@>\rho'>>(K\GD, V_{K\GD}(a))\\
@V(\R_+s\times \id)VV @VsVV\\
(\R_+\GD(a)\times U_K^n,\, \rho^{-1}(V_{K\GD(a)}(a)))
@>\rho>>(K\GD(a), V_{K\GD(a)}(a))
\end{CD}
$$
appears. Here $s$ is the natural regular map resolving
singularities of $K\GD(a)$, $\rho$ and $\rho'$ are natural
projections and $\R_+s$ is a map $\R_+\GD\to\R_+\GD(a)$ defined
by $s$. The preimage of $K\GG$ under $\rho$ is
the union of
$K\GS$ with $\GS\in\mathcal{G}'(\GD)$ and $C_\GD(\GS)\supset
C_{\GD(a)}(\GG)$. By \ref{2.5.D}, the set of singular points
of $V_{K\GD}(a)$
does not intersect $K\GS$, and $V_{K\GD}(a)$ is transversal to
$K\GS$.

If $\GS\in\mathcal{G}'(\GD)$, $C_\GD(\GS)\supset C_{\GD(a)}(\GG)$ and
$\dim \GS=\dim \GG$, then  $\R_+s$ defines an isomorphism
$\R_+C_\GS(\GS)\to\R_+C_\GG(\GG)$, and if
$\GS\in\mathcal{G}'(\GD)$, $C_\GD(\GS)\supset C_{\GD(a)}(\GG)$ and
$\dim \GS>\dim \GG$, then $\R_+s$ defines a map
$\R_+C_\GS(\GS)\to\R_+C_\GG(\GG)$ which is a factorization by the
action of quasi-homotheties $qh_{w,t}$ with $w\in
C_\GS(\GS)$, $w\perp C_\GG(\GG)$. By \ref{2.5.E}, in the latter case
variety $V_{K\GS}(a^\GG)$ coinciding, by \ref{2.5.A}, with
$V_{K\GD}(a)\cap K\GS$ is invariant under the same
quasi-homotheties. Hence $V_{K\GD}(a)=s^{-1}V_{K\GD(a)}(a)$ and
hypersurface $\rho^{-1}V_{K\GD}(a)$, being the image of
$\rho^{\prime-1}V_{K\GD}(a)$ under $\R_+\times \id$, appears to be
nonsingular along its intersection with  $\R_+\GG\times U_K^n$ and
transversal to $\R_+\GG\times U_K^n$.
\end{proof}

%
%



\section{Patchworking}\label{s4}

\subsection{Patchworking  L-polynomials}\label{s4.1}Let
$\GD$, $\GD_1$,
\dots, $\GD_s\subset\R^n$ be (convex integer) polyhedra with
$\GD=\bigcup_{i=1}^s\GD_i$ and $\Int\GD_i\cap \Int\GD_j=\empt$ for
$i\ne j$. Let $\nu:\GD\to\R$ be a nonnegative convex function
satisfying to the following conditions:
\begin{enumerate}
\item\label{1nu} all the restrictions $\nu|_{\GD_i}$ are linear;
\item\label{2nu} if the restriction of $\nu$ to an open set is linear
then this set is contained in one of $\GD_i$;
\item\label{3nu} $\nu(\GD\cap\Z^n)\subset\Z$.
\end{enumerate}
\begin{rem}\label{}Existence of such a function $\nu$ is a restriction on
a collection $\GD_1,\dots,\GD_s$. For example, the collection of convex
polygons shown in Figure \ref{f24} does not admit such a
function.\end{rem}
\begin{figure}[h]
\centerline{\epsffile{pw-f24.eps}}
\caption{}
\label{f24}
\end{figure}

Let $a_1,\dots,a_s$ be Laurent polynomials over $K$ in $n$ variables
with $\GD(a_i)=\GD$. Let $a_i^{\GD_i\cap\GD_j}=a_j^{\GD_i\cap\GD_j}$
for any $i$, $j$. Then, obviously, there exists an unique L-polynomial $a$
with $\GD(a)=\GD$ and $a^{\GD_i}=a_i$ for $i=1,\dots,s$. If
$a(x_1,\dots,x_n)=\sum_{\Go\in\Z^n}a_\Go x^\Go$, we put
$b(x,t)=\sum_{\Go\in\Z^n}a_\Go x^\Go t^{\nu(\Go)}$.
This L-polynomial in $n+1$ variables is considered below also as a
one-parameter family of L-polynomials in $n$ variables. Therefore
let me introduce the corresponding notation: put
$b_t(x_1,\dots,x_n)=b(x_1,\dots,x_n,t)$.
L-polynomials $b_t$ are said to be obtained by {\it patchworking}
L-polynomials $a_1,\dots,a_s$ by $\nu$ or, briefly, $b_t$ is a
patchwork of L-polynomials $a_1,\dots,a_s$ by $\nu$.

\subsection{Patchworking  charts}\label{s4.2}Let $a_1,\dots,a_s$ be
Laurent polynomials over $K$ in $n$ variables with $\Int\GD(a_i)\cap
\Int\GD(a_j)=\empt$ for $i\ne j$. A pair $(\GD\times U_K^n,\, \upsilon)$
is said to be  obtained by {\it patchworking} $K$-charts
of Laurent polynomials $a_1,\dots,a_s$ and it is a {\it patchwork} of
$K$-charts
of L-polynomials $a_1,\dots,a_s$ if $\GD=\bigcup_{i=1}^s\GD(a_i)$ and
one can choose $K$-charts $(\GD(a_i)\times U_K^n, \upsilon_i)$ of Laurent
polynomials $a_1,\dots,a_s$ such that
$\upsilon=\bigcup_{i=1}^s\upsilon_i$.

\subsection{The Main Patchwork Theorem}\label{s4.3}Let $\GD$, $\GD_1$,
\dots, $\GD_s$, $\nu$, $a_1$, \dots, $a_s$, $b$ and $b_t$ be as in
Section \ref{s4.1} ($b_t$ is a patchwork of L-polynomials
$a_1$, \dots, $a_s$ by $\nu$).

\begin{prop}\label{4.3.A}If L-polynomials $a_1,\dots,a_s$ are completely
nondegenerate then there exists $t_0>0$ such that for any $t\in(0,t_0]$
a $K$-chart of L-polynomial $b_t$ is obtained by patchworking
$K$-charts of L-polynomials $a_1,\dots,a_s$.
\end{prop}
\begin{proof}Denote by $\mathcal G$ the union $\bigcup_{i=1}^s\mathcal{G}(\GD_i)$. For $\GG\in\mathcal G$ denote by $\tilde\GG$ the graph of
$\nu|_\GG$. It is clear that $\GD(b)$ is the convex hull of graph
of $\nu$, so $\tilde\GG\in\mathcal{G}(\GD(b))$ and thus there is an
injection $\mathcal{G}\to\mathcal{G}(\GD(b)):\GG\mapsto\tilde\GG$.
Restrictions $\tilde\GG\to\GG$ of the natural projection $pr:\R^{n+1}\to\R^n$ are homeomorphisms, they are denoted by $g$.

Let $p:\GD(b)\times U_K^{n+1}\to K\GD(b)$ be the composition of
the homeomorphism
$$\begin{CD}\GD(b)\times U_K^{n+1}@>h\times \id>>\R_+\GD(b)\times U_K^{n+1}
\end{CD}$$
and the natural projection $\rho:\R_+\GD(b)\times U_K^{n+1}\to K\GD(b)$
(cf. Section \ref{s3.3}), so
the pair $\left(\GD(b)\times U_K^{n+1},\, p^{-1}V_{K\GD(b)}(b)\right)$
is a $K$-chart of $b$. By \ref{2.5.A}, for $i=1,\dots,s$ the pair
$$\left(\widetilde{\GD(a_i)}\times U_K^{n+1},
p^{-1}(V_{K\GD(b)}(b)\cap\widetilde{\GD(a_i)}\times U_K^n)\right)$$
is a $K$-chart of L-polynomial $b^{\widetilde{\GD(a_i)}}$.

The pair
$$\left(\widetilde{\GD(a_i)}\times U_K^n,\,
p^{-1}(V_{K\GD(b)}(b)\cap\widetilde{\GD(a_i)}\times U_K^n)\right)$$
which is cut out by this pair
on $\widetilde{\GD(a_i)}\times U_K^n$   is transformed
 by $g\times \id:\widetilde{\GD(a_i)}\times
U_K^n\to\GD(a_i)\times U_K^n$ to a $K$-chart of $a_i$. Indeed,
$g:\widetilde{\GD(a_i)}\to\GD(a_i)$ defines an isomorphism
$g^*:K\GD(a_i)\to K\widetilde{\GD(a_i)}$ and since
$b^{\widetilde{\GD(a_i)}}(x_1,\dots,x_n,1)=a_i(x_1,\dots,x_n)$,
it follows that
$g^*:V_{K\GD(a_i)}(a_i)=V_{K\widetilde{\GD(a_i)}}
(b^{\widetilde{\GD(a_i)}})$ and  $g$ defines a homeomorphism
of the pair
$\left(\widetilde{\GD(a_i)}\times U_K^n,\,
p^{-1}(V_{K\GD(b)}(b)\cap\widetilde{\GD(a_i)}\times U_K^n)\right)$ to a
$K$-chart of L-polynomial $a_i$.

Therefore the pair
$$\left(\bigcup_{i=1}^s\widetilde{\GD(a_i)}\times U_K^n,\,
p^{-1}(V_{K\GD(b)}(b)\cap\bigcup_{i=1}^s\widetilde{\GD(a_i)}\times
U_K^n)\right)$$
is a result of patchworking $K$-charts of $a_1,\dots,a_s$.

For $t>0$ and $\GG\in\mathcal{G}'(\GD)$ let us construct a ring
homomorphism
$$K[C_{\GD(b)}(\GD(b)\cap pr^{-1}(\GG))\cap\Bbb
Z^{n+1}]\to K[C_\GD(\GG)\cap\Z^n]$$
which maps a monomial $x_1^{\Go_1}\dots x_n^{\Go_n}x_{n+1}^{\Go_{n+1}}$
to $t^{\Go_{n+1}}x_1^{\Go_1}\dots x_n^{\Go_n}$. This homomorphism
corresponds to the embedding
$$KC_\GD(\GG)\to KC_{\GD(b)}(\GD(b)\cap pr^{-1}(\GG))$$
extending the embedding
$$K\R^n\to K\R^{n+1}:(x_1,\dots,x_n)\mapsto(x_1,\dots,x_n,t)$$.
The embeddings constructed in this way agree  to each other and define
an embedding $K\GD\to K\GD(b)$. Denote the latter embedding by $i_t$.
It is clear that $V_{K\GD}(b_t)=i_t^{-1}V_{K\GD(b)}(b)$.

The sets $\rho^{-1}i_tK\GD$ are smooth hypersurfaces of
$\GD(b)\times U_K^{n+1}$, comprising a smooth isotopy. When $t\to0$,
the hypersurface $\rho^{-1}i_tK\GD$ tends (in $C^1$-sense) to
$$\bigcup_{i=1}^s\widetilde{\GD(a_i)}\times U_K^n.$$
By \ref{3.3.A}, $\rho^{-1}V_{K\GD(b)}(b)$ is transversal to
each of $$\R_+\widetilde{\GD(a_i)}\times U_K^{n+1}$$
and hence, the intersection
$\rho^{-1}(i_tK\GD)\cap\rho^{-1}(V_{K\GD(b)}(b))$ for sufficiently
small $t$  is mapped to
$$V_{K\GD(b)}(b)\cap\bigcup_{i=1}^s\widetilde{\GD(a_i)}\times
U_K^n$$
by some homeomorphism
$$\rho^{-1}i_tK\GD\to\bigcup_{i=1}^s\widetilde{\GD(a_i)}\times U_K^n.$$
Thus the pair
$$\left(\rho^{-1}i_tK\GD,\,
\rho^{-1}i_tK\GD\cap\rho^{-1}V_{K\GD(b)}(b)\right)$$
is a result of patchworking $K$-charts of L-polynomials $a_1,\dots,a_s$
if $t$ belongs to a segment of the form $(0,t_0]$. On the other hand,
since $V_{K\GD}(b_t)=i^{-1}_tV_{K\tilde\GD}(b)$,
$$\rho^{-1}i_tK\GD\cap\rho^{-1}V_{K\GD(b)}(b)=\rho^{-1}i_tV_{K\GD}(b_t)$$
and, hence, the pair
$$\left(\rho^{-1}i_tK\GD,\,\rho^{-1}i_tK\GD\cap\rho^{-1}V_{K\GD(b)}(b)\right)$$
is homeomorphic to a $K$-chart of L-polynomial $b_t$.
\end{proof}


\def\up{\upsilon}

\section{Perturbations smoothing a singularity of
hypersurface}\label{s5}

The construction of the previous Section can be interpreted as a purposeful
smoothing of an algebraic hypersurface with singularities, which results
in replacing
of neighborhoods of singular points by new fragments of
hypersurface, having a prescribed topological structure (cf. Section
\ref{s1.9}). According to well known theorems of theory of
singularities, all theorems on singularities of algebraic hypersurfaces
are extended to singularities of significantly wider class of
hypersurfaces. In particular, the construction of perturbation based on
patchworking is applicable in more general situation. For singularities
of simplest types this construction  together with some results of
topology of algebraic curves allows to get a topological classification
of perturbations  which smooth singularities completely.

The aim if this Section is to adapt patchworking to needs of singularity
theory.

\subsection{Singularities of hypersurfaces}\label{s5.1}Let
$G\subset K^n$ be an
open set, and let $\varphi:G\to K$ be an analytic function. For
$U\subset G$  denote by $V_U(\varphi)$  the set
$\{x\in U\,|\,\varphi(x)=0\}$.

By  {\it singularity} of a hypersurface $V_G(\varphi)$ at the point
$x_0\in V_G(\varphi)$ we mean the class of germs of hypersurfaces which
are diffeomorphic to the germ of  $V_G(\varphi)$ at $x_0$.
In other words,  hypersurfaces $V_G(\varphi)$ and $V_H(\psi)$ have the
same singularity at  points $x_0$ and $y_0$, if there exist
neighborhoods $M$ and $N$ of $x_0$ and $y_0$ such
that the pairs $(M, V_M(\varphi))$, $(N, V_N(\psi))$ are diffeomorphic.
When considering a singularity of hypersurface at a point $x_0$, to
simplify the formulas we shall assume that $x_0=0$.

The {\it multiplicity} or the {\it Milnor number} of a hypersurface
$V_G(\varphi)$ at $0$ is the dimension
$$dim_KK[[x_1,\dots,x_n]]/(\p f/\p x_1,\dots,
\p f/\p x_n)$$
of the quotient of the formal power series ring by the ideal generated
by  partial derivatives $\p f/\p x_1,\dots,\p
f/\p x_n$ of the Taylor series expansion $f$ of the function $\varphi$ at $0$.
This number is an invariant of the singularity (see \cite{3s}). If it is
finite, then we say that the singularity is of {\it finite
multiplicity}.

If the singularity of $V_G(\Gf)$ at $x_0$ is of finite
multiplicity, then this singularity is isolated, i.e.  there exists a
neighborhood $U\subset K^n$ of $x_0$, which does not contain singular
points of $V_G(\varphi)$.  If $K=\CC$ then the converse is true: each
isolated singularity of a hypersurface is of finite multiplicity.  In
the case of isolated singularity, the boundary of a ball
$B\subset K^n$, centered at $x_0$ and small enough, intersects
$V_G(\varphi)$ only at nonsingular points and only transversely, and
 the pair $(B, V_B(\varphi))$ is homeomorphic to the cone over its
boundary $(\p B, V_{\p B}(\varphi))$ (see \cite{15s}, Theorem 2.10).
In such a case the pair $(\p B, V_{\p B}(\varphi))$ is called the {\it
link} of singularity of $V_G(\varphi)$ at $x_0$.

The following Theorem shows that the class of singularities of finite
multiplicity of analytic hypersurfaces coincides with the class of
singularities of finite multiplicity of algebraic hypersurfaces.

\begin{prop}[Tougeron's theorem]\label{5.1.A}(\em{see, for example,}
\cite{3s}, Section \ref{s6.3}). If the singularity at $x_0$ of a
hypersurface $V_G(\varphi)$ has finite Milnor number $\mu$, then there
exist a neighborhood $U$ of\, $x_0$ in $K^n$ and a diffeomorphism $h$
of this neighborhood onto a neighborhood of\,  $x_0$ in $K^n$ such that
$h(V_U(\varphi))=V_{h(U)}(f_{(\mu+1)})$, where $f_{(\mu+1)}$
is the Taylor polynomial of $\varphi$ of degree $\mu+1$ .
\end{prop}

The notion of Newton polyhedron is extended  over in a natural way to power
series. The Newton polyhedron $\GD(f)$ of the series
$f(x)=\sum_{\Go\in\Z^n}a_\Go x^\Go$ (where $x^\Go=
x_1^{\Go_1}x_2^{\Go_2}\dots x_n^{\Go_n}$) is the convex hull
of the set $\{\Go\in\R^n\,|\,a_\Go\ne0\}$. (Contrary to
 the case of a
polynomial, the Newton polyhedron $\GD(f)$ of a power series may have
infinitely many faces.)

However in the singularity theory the notion
of  Newton diagram occurred to be more important. The {\it Newton
diagram} $\GG(f)$ of a power series $f$ is the union of the proper faces
of the Newton polyhedron which face the origin, i.e. the union of the
faces $\GG\in\mathcal{G}'(\GD(f))$ for which cones $D^+C_{\GD(f)}(\GG)$
intersect the open positive orthant $\Int A^n=\{x\in\R^n\,|\,x_1>0,
\dots, x_n>0\}$.

It follows from the definition of the Milnor number  that, if the
singularity of $V_G(\varphi)$ at $0$ is of finite multiplicity, the
Newton diagram of the Taylor series of $\varphi$ is compact, and its
distance from each of the coordinate axes is at most $1$.

For a power series $f(x)=\sum_{\Go\in\Z^n}f_\Go x^\Go$ and
a set $\GG\subset \R^n$ the power series $\sum_{\Go\in\GG\cap\Bbb
Z^n}f_\Go x^\Go$ is called  $\GG$-{\it truncation} of $f$ and
denoted by $f^\GG$ (cf. Section \ref{s2.1}).

Let the  Newton diagram of the Taylor series $f$ of a function
$\varphi$ be compact. Then $f^{\GG(f)}$  is a polynomial. The pair
$(\GG(f)\times U_K^n, \Gg)$ is said to be a {\it nonreduced chart}
of germ of hypersurface $V_G(\varphi)$ at $0$ if there exists  a
$K$-chart $(\GD(f^{\GG(f)}\times U_K^n, \up)$ of  $f^{\GG(f)}$ such
that $\gamma=\up\cap(\GG(f)\times U_K^n)$. It is clear that a nonreduced
chart of germ of hypersurface is comprised of $K$-charts of $f^\GG$,
where $\GG$ runs over the set of all faces of the Newton diagram.

A power series $f$ is said to be {\it nondegenerate} if its Newton
diagram is compact and the distance between it and each of the
coordinate axes is at most $1$ and for any its face $\GG$ polynomial
$f^\GG$ is completely nondegenerate. In this case about the germ of
$V_G(\varphi)$ at zero we
say that {\it it is placed nondegenerately}. It is not difficult
to prove that nondegenerately placed germ defines a singularity of
 finite multiplicity. It is convenient to place the charts of germs
of hypersurfaces  in $K^n$ by
a natural map $A^n\times U_K^n\to K^n:(x,y)\mapsto S_y(x)$
(like $K$-charts of an L-polynomial, cf. Section \ref{s3.3}).
Denote by $\Sigma_K(\varphi)$ the image of $\GG(f)\times
U_K^n$ under this map; the image of nonreduced chart of germ of
hypersurface $V_G(\varphi)$ at zero under this map is called a
{\it (reduced) chart of germ\/} of $V_G(\varphi)$ at the origin. It
follows from Tougeron's theorem that in this case adding a monomial of
the form $x_i^{m_i}$ to $\varphi$ with $m_i$ large enough does not
change the singularity. Thus, without changing the singularity, one can
make the Newton diagram meeting the coordinate axes.

In the case when this takes place and the Taylor series of $\varphi$ is
nondegenerate there exists a ball $U\subset K^n$ centered at $0$
such that the pair $(U, V_U(\varphi))$ is homeomorphic to the cone over a
chart of germ of $V_G(\varphi)$. This follows from Theorem \ref{5.1.A} and
from results of Section \ref{s2.5}.

Thus if the Newton diagram meets all coordinate axes and the Taylor series
of $\varphi$ is nondegenerate, then the chart of germ of $V_G(\varphi)$
at zero is homeomorphic to the link of the singularity.

\subsection{Evolving of a singularity}\label{s5.2}Now let the function
$\Gf:G\to K$ be included as $\Gf_0$ in a family of analytic functions
$\Gf_t:G\to K$ with $t\in[0,t_0]$. Suppose that this is an analytic
family in the sense that the function
$G\times[0,t_0]\to K:(x,t)\mapsto\Gf_t(x)$ which is  determined by it
is real analytic. If the
hypersurface $V_G(\Gf)$ has an isolated singularity at $x_0$, and if
there exists a neighborhood $U$ of $x_0$ such that the hypersurfaces
$V_G(\Gf_t)$ with $t\in[0,t_0]$  have no singular points in $U$,
then the family of functions $\Gf_t$ with $t\in[0,t_0]$
is said to {\it evolve} the singularity of $V_G(\Gf)$ at $x_0$.

If the family $\Gf_t$ with $t\in[0,t_0]$ evolves the singularity of
the hypersurface $V_G(\Gf_0)$ at $x_0$, then there exists a ball
$B\subset K^n$ centered at $x_0$ such that
\begin{enumerate}
\item\label{ev1} for $t\in[0,t_0]$ the sphere $\p B$ intersects
$V_G(\Gf_t)$ only in nonsingular points of the hypersurface and
only transversely,
\item\label{ev2} for $t\in(0,t_0]$ the ball $B$ contains no singular
point of the hypersurface $V_G(\Gf_t)$,
\item\label{ev3} the pair $(B, V_B(\Gf_0))$ is homeomorphic to the cone
over its boundary $(\p B$, $V_{\p B}(\Gf_0))$.
\end{enumerate}

Then the family of pairs $(B, V_B(\Gf_t))$ with $t\in[0,t_0]$ is called
an {\it evolving} of the germ of $V_G(\Gf_0)$ in
$x_0$. (Following the standard terminology of the
singularity theory, it would be more correct to say not a on family of
pairs, but rather a family of germs or even germs of a family;
however, from the topological viewpoint, which is more natural in
the context of the topology of real algebraic varieties, the distinction
between a family of pairs satisfying \ref{ev1} and \ref{ev2} and
the corresponding family of germs is of no importance, and so we shall
ignore it.)

Conditions \ref{ev1} and \ref{ev2} imply  existence of a
smooth isotopy $h_t:B\to B$ with $t\in(0,t_0]$, such that $h_{t_0}=\id$
and $h_t(V_B(\Gf_{t_0}))=V_B(\Gf_t)$, so that the pairs
$(B,V_B(\Gf_t))$ with $t\in(0,t_0]$ are homeomorphic to each other.

Given germs determining the same singularity, a evolving
of one of them obviously corresponds to a diffeomorphic evolving of
the other germ. Thus, one  may speak not only of
evolvings of germs, but also of {\it evolvings of singularities}
of a hypersurface.

The following three topological classification questions
on evolvings arise.

\begin{prop}\label{5.2.A}Up to homeomorphism, what manifolds  can appear as
$V_B(\Gf_t)$ in evolvings of a given singularity?
\end{prop}

\begin{prop}\label{5.2.B}Up to homeomorphism, what pairs can appear as
$(B,V_B(\Gf_t))$ in evolvings of a given singularity?
\end{prop}

Smoothings $(B,V_B(\Gf_t))$ with $t\in[0,t_0]$ and
$(B',V_{B'}(\Gf'_t))$ with $t\in[0,t'_0]$ are said to be {\it
topologically equivalent} if there exists an isotopy $h_t:B\to B'$ with
$t\in[0, \min(t_0,t'_0)]$, such that $h_0$ is a diffeomorphism and
$V_{B'}(\Gf'_t)=h_tV_B(\Gf_t)$ for $t\in[0, \min(t_0,t'_0)]$.

\begin{prop}\label{5.2.C}Up to topological equivalence, what are the evolvings
of a given singularity?
\end{prop}

Obviously, \ref{5.2.B} is a refinement of \ref{5.2.A}. In turn,
\ref{5.2.C} is more refined than \ref{5.2.B}, since in \ref{5.2.C} we are
interested not only in the type of the pair obtained in result of the
evolving, but also the manner in which the pair is attached to the
link of the singularity.

In the case $K=\R$ these questions have been answered
in literature only for several simplest
singularities.

In the case $K=\C$ a evolving of a given singularity is unique
from each of the  three points of view, and there is an extensive literature
(see, for example, \cite{18d}) devoted to its topology (i.e., questions
\ref{5.2.A} and \ref{5.2.B}).

By the way, if we want to get questions for $K=\C$ which are
truly analogous to questions \ref{5.2.A} --- \ref{5.2.C} for $K=\R$, then we
have to replace evolvings by deformations with nonsingular fibers and
one-dimensional complex bases, and the variety $V_B(\Gf_t)$ and the
pairs $(B, V_B(\Gf_t))$ have to be considered along with the monodromy
transformations. It is reasonable to suppose that there are interesting
connections between questions \ref{5.2.A} --- \ref{5.2.C} for a real
singularity and their counter-parts  for the complexification of the
singularity.

\subsection{Charts of evolving}\label{s5.3}Let the Taylor series $f$ of
function $\Gf:G\to K$ be nondegenerate and its Newton diagram meets
all the coordinate axes. Let  a family of functions $\Gf_t:G\to K$ with
$t\in[0,t_0]$ evolves the singularity of $V_G(\Gf)$ at $0$. Let $(B,
V_B(\Gf_t))$ be the corresponding evolving of the germ of this
hypersurface and $h_t:B\to B$ with $t\in(0,t_0]$ be an isotopy with
$h_{t_0}=\id$ and $h_t(V_B(\Gf_{t_0}))=V_B(\Gf_t))$ existing by
conditions \ref{ev1} and \ref{ev2} of the previous Section.
Let $(\Sigma_K(\Gf),\gamma)$ be a chart of germ of hypersurface
$V_G(\Gf)$ at zero and $g:(\Sigma_K(\Gf),\gamma)\to(\p B,
V_{\p B}(\Gf))$ be the natural homeomorphism of it to link of
the singularity.

Denote by $\Pi_K(\Gf)$ a part of $K^n$ bounded by $\Sigma_K(\Gf)$. It
can be presented as a cone over $\Sigma_K(\Gf)$ with vertex at zero.

One can choose the isotopy $h_t:B\to B$, $t\in(0,t_0]$
such that its restriction to $\p B$ can be extended to an isotopy
$h'_t:\p B\to\p B$ with $t\in[0,t_0]$ (i.e., extended for $t=0$).

We shall call the pair $(\Pi_K(\Gf)$, $\tau)$ a {\it chart of
evolving} $(B$, $V_B(\Gf_t))$, $t\in[0,t_0]$, if there exists a
homeomorphism $(\Pi_K(\Gf)$, $\tau)\to(B$, $V_B(t_0))$, whose
restriction $\Sigma_K(\Gf)\to\p B$ is the composition
$\begin{CD}\Sigma_K(\Gf)@>g>>\p B@>h'_0>>\p B\end{CD}$. One can see
that the boundary $(\p\Pi_K(\Gf)$, $\p\tau)$ of
a chart of evolving is a chart $(\Sigma_K(\Gf)$, $\gamma)$ of the germ
of the hypersurface at zero, and a chart of evolving is a pair
obtained by evolving which is glued to $(\Sigma_K(\Gf)$, $\gamma)$ in
natural way. Thus that the chart of an evolving
describes the evolving up to topological equivalence.

\subsection{Construction of evolvings by patchworking}\label{s5.4}Let
the Taylor series $f$ of function $\Gf:G\to K$ be nondegenerate and its
Newton diagram $\GG(f)$ meets all the coordinate axes.

Let $a_1, \dots, a_s$ be completely nondegenerate polynomials over
$K$ in $n$ variables with $\Int\GD(a_i)\cap \Int\GD(a_j)=\varnothing$
and $a_i^{\GD(a_i)\cap\GD(a_j)}=a_j^{\GD(a_i)\cap\GD(a_j)}$ for $i\ne
j$.  Let $\bigcup_{i=1}^s\GD(a_i)$ be the polyhedron bounded
by the coordinate axes and Newton diagram $\GG(f)$. Let
$a_i^{\GD(a_i)\cap\GD(f)}=f^{\GD(a_i)\cap\GD(f)}$ for $i=1,\dots, s$.
Let $\nu:\bigcup_{i=1}^s\GD(a_i)\to\R$ be a nonnegative convex
function which is equal to zero on $\GG(f)$ and satisfies conditions
\ref{1nu}, \ref{2nu}, \ref{3nu} of Section \ref{s4.1}
with polyhedra $\GD(a_1)$, \dots, $\GD(a_s)$. Then polynomials
$a_1, \dots, a_s$ can be "glued to $\Gf$ by $\nu$" in the following
way generalizing patchworking L-polynomials of Section \ref{s4.1}.

Denote by $a$ the polynomial defined by conditions $a^{\GD(a_i)}=a_i$ for
$i=1, \dots, s$ and $a^{\bigcup_{i=1}^{s}\GD(a_i)}=a$. If
$a(x)=\sum_{\Go\in\Z^n}a_\Go x^\Go$ then we put
$$\Gf_t(x)=\Gf(x)+(\sum_{\Go\in\Z^n}a_\Go x^\Go
t^{\nu(\Go)})-a^{\GG(f)}x.$$

\begin{prop}\label{5.4.A} Under the conditions above there exists $t_0>0$ such that
the family of functions $\Gf_t:G\to K$ with $t\in[0,t_0]$ evolves the
singularity of $V_G(\Gf)$ at zero. The chart of this evolving is
patchworked from $K$-charts of $a_1, \dots, a_s$.  \end{prop}

In the case when $\Gf$ is a polynomial, Theorem \ref{5.4.A} is a slight
modification of a special case of Theorem \ref{4.3.A}. Proof of \ref{4.3.A}
is easy to transform to the proof of this version of \ref{5.4.A}. The general
case can be reduced to it by
Tougeron Theorem, or one can prove it directly, following to scheme of
proof of Theorem \ref{4.3.A}.\qed

We shall call the evolvings obtained by the scheme described in this
Section  {\it patchwork} evolvings.

\def\Gf{\varphi}
\def\neis{neighborhoods }
\def\nei{neighborhood }
\def\GD{\Delta}
\def\GG{\Gamma}
\def\Ge{\varepsilon}

\section{Approximation of hypersurfaces of $K\R^n$}\label{s6}

\subsection{Sufficient truncations}\label{s6.1}Let $M$ be a smooth
submanifold of a smooth manifold $X$. Remind that by a {\it tubular
neighborhood} of $M$ in $X$ one calls a submanifold $N$ of $X$ with
$M\subset \Int N$ equipped  with a {\it tubular fibration}, which  is a
smooth retraction $p:N\to M$ such that for any point $x\in M$ the
preimage $p^{-1}(x)$ is a smooth submanifold of $X$
diffeomorphic to $D^{\dim X-\dim M}$. If $X$ is equipped  with a
metric and each fiber of the tubular fibration $p:N\to M$ is contained
in a ball of radius $\Ge$ centered in the point of
intersection of the fiber with $M$, then $N$ is called  a {\it tubular
$\Ge$-neighborhood} of $M$ in $X$.

We need tubular neighborhoods mainly for formalizing a notion of
approximation of a submanifold by a submanifold. A manifold presented
as the image of a smooth section of the tubular fibration of a tubular
$\Ge$-neighborhood of $M$ can be considered as sufficiently
close  to $M$:
it is naturally isotopic to $M$ by an isotopy moving each point at
most by $\Ge$.

We shall consider the space $\R^n\times U_K^n$ as a flat Riemannian
manifold with metric defined by the standard  Euclidian metric of
$\R^n$ in the case of $K=\R$ and by the
standard Euclidian metric of $\R^n$ and the standard flat metric of
the torus $U_{\C}^n=(S^1)^n$  in the case of $K=\C$.

An $\Ge$-sufficiency of truncations of Laurent polynomial defined below
and the whole theory related with this notion presuppose that it has
been chosen a class of tubular neighborhoods of smooth submanifolds of
$\R^n\times U_K^n$ invariant under translations
$T_\Go\times \id_{U_K^n}$ and that for any two
tubular \neis $N$ and $N'$ of the same $M$, which belong to this class,
restrictions of tubular fibrations $p:N\to M$ and $p':N'\to M$ to
$N\cap N'$ coincide. One of such classes is the collection of all {\it
normal} tubular neighborhoods, i.e. tubular \neis with fibers
consisting of segments of geodesics which start from the same point of
the submanifold in directions orthogonal to the submanifold. Another
class, to which we shall turn in Sections \ref{s6.7} and \ref{s6.8}, is
the class of tubular \neis whose fibers lie in  fibers $\R^{n-1}\times
t\times U_K^{n-1}\times s$ of $\R^n\times U_K^n$ and consist of
segments of geodesics which are orthogonal to intersections of the
corresponding manifolds with these $\R^{n-1}\times t\times
U_K^{n-1}\times s$. The intersection of such a tubular \nei of $M$ with
the fiber $\R^{n-1}\times t\times U_K^{n-1}\times s$ is a normal
tubular \nei of $M\cap(\R^{n-1}\times t\times U_K^{n-1}\times s)$ in
$\Bbb R^{n-1}\times t\times U_K^{n-1}\times s$. Of course, only
manifolds transversal to $\R^{n-1}\times t\times U_K^{n-1}\times s$
have tubular \neis of this type.

Introduce a norm in vector space of Laurent polynomials over $K$ on $n$
variables:
$$||\sum_{\Go\in\Z^n}a_\Go
x^\Go||=\max\{|a_\Go|\,|\,\Go\in\Z^n\}.$$

Let $\GG$ be a subset of $\R^n$ and $\Ge$  a positive number. Let
$a$ be a Laurent polynomial over $K$ in $n$ variables and $U$  a
subset of $K\R^n$. We shall say that {\it in $U$ the truncation
$a^\GG$ is $\Ge$-sufficient for} $a$ (with respect to the chosen class of
tubular neighborhoods), if for any Laurent polynomial $b$ over $K$
satisfying the conditions $\GD(b)\subset\GD(a)$, $b^\GG=a^\GG$ and
$||b-b^\GG||\le||a-a^\GG||$ (in particular, for $b=a$ and $b=a^\GG$) the
following condition takes place: \begin{enumerate}
\item $U\cap SV_{K\R^n}(b)=\varnothing$,
\item the set $la(U\cap V_{K\R^n}(b))$ lies in a tubular $\Ge$-\nei $N$
(from the chosen class) of $la(V_{K\R^n}(a^\GG)\sminus SV_{K\R^n}(a^\GG))$
and
\item $la(U\cap V_{K\R^n}(b))$ can be extended to the image of a smooth
section of the tubular fibration $N\to la(V_{K\R^n}(a^\GG)\smallsetminus
SV_{K\R^n}(a^\GG))$.
\end{enumerate}

The $\Ge$-sufficiency of $\GG$-truncation of Laurent polynomial $a$ in
$U$ means, roughly speaking, that monomials which are not in $a^\GG$
have a small influence on $V_{K\R^n}(a)\cap U$.

\begin{prop}\label{6.1.A}If $a^\GG$ is $\Ge$-sufficient for $a$ in open
sets $U_i$
with $i\in\mathcal J$, then it is $\Ge$-sufficient for $a$ in
$\bigcup_{i\in\mathcal{J}}U_i$ too.\qed
\end{prop}

Standard arguments based on Implicit Function Theorem give the
following Theorem.

\begin{prop}\label{6.1.B}If a set $U\subset K\R^n$ is compact and contains no
singular points of a hypersurface $V_{K\R^n}(a)$, then for any
tubular \nei $N$ of $V_{K\R^n}(a)\sminus SV_{K\R^n}(a)$ and any
polyhedron $\GD\supset\GD(a)$
there exists  $\Gd>0$ such that for any Laurent
polynomial $b$ with $\GD(b)\subset\GD$ and $||b-a||<\Gd$
the hypersurface $V_{K\R^n}(b)$ has no singularities in
$U$, intersection $U\cap V_{K\R^n}(b)$ is contained in
$N$ and can be extended to the image of a smooth section of a
tubular fibration $N\to V_{K\R^n}(a)\sminus SV_{K\R^n}(a)$.
\qed\end{prop}

From this the following proposition follows easily.

\begin{prop}\label{6.1.C}If $U\in K\R^n$ is compact and
$a^\GG$ is $\Ge$-sufficient truncation of $a$ in $U$, then for any
polyhedron $\GD\supset\GD(a)$ there exists $\Gd>0$ such that for any
Laurent polynomial $b$ with $\GD(b)\subset\GD$,
$b^\GG=a^\GG$ and $||b-a||<\Gd$ the truncation $b^\GG$ is
$\Ge$-sufficient in $U$.  \qed\end{prop}

In the case of $\GG=\GD(a)$ proposition \ref{6.1.C} turns to the following
proposition.

\begin{prop}\label{6.1.D}If a set $U\subset K\R^n$ is compact and contains no
singular points of $V_{K\R^n}(a)$ and $la (V_{K\R^n}(a))$ has a tubular
\nei of the chosen type, then for
any $\Ge>0$ and any polyhedron $\GD\supset\GD(a)$ there exists
$\Gd>0$ such that for any Laurent polynomial $b$ with
$\GD(b)\subset\GD$, $||b-a||<\Gd$ and $b^{\GD(a)}=a$ the
truncation $b^{\GD(a)}$ is $\Ge$-sufficient in $U$.
\qed\end{prop}

The following proposition describes behavior of the $\Ge$-sufficiency
under quasi-homotheties.

\begin{prop}\label{6.1.E}Let $a$ be a Laurent polynomial over $K$ in $n$
variables. Let $U\subset K\R^n$, $\GG\subset\R^n$, $w\in\R^n$.
Let $\Ge$ and $t$ be positive numbers. Then $\Ge$-sufficiency of
$\GG$-truncation $a^\GG$ of $a$ in $qh_{w,t}(U)$ is equivalent to
$\Ge$-sufficiency of $\GG$-truncation of $a\circ qh_{w,t}$ in $U$.
\end{prop}

The proof follows from comparison of the definition of $\Ge$-sufficiency
and the following two facts. First, it is obvious that
$$qh_{w,t}(U)\cap V_{K\R^n}(b)=qh_{w,t}(U\cap
qh^{-1}_{w,t}(V_{K\R^n}(b)))=qh_{w,t}(U\cap
V_{K\R^n}(b\circ qh_{w,t})),$$
and second, the transformation $T_{(\ln t)w}\times \id_{U_K^n}$ of
$\R^n\times U_K^n$ corresponding, by \ref{2.1.A}, to $qh_{w,t}$ preserves
the chosen  class of
tubular
$\Ge$\neis.\qed

\subsection{Domains of $\Ge$-sufficiency of face-truncation}\label{s6.2}
For $A\subset\R^n$ and $B\subset K\R^n$  denote by $qh_A(B)$ the union
$\bigcup_{\Go\in A}qh_\Go(B)$.

For $A\subset\R^n$ and $\rho>0$ denote by $\frak N_\rho(A)$
the set $\{x\in\R^n|dist(x,A)<\rho\}$.

For $A,B\subset\R^n$ and $\lambda\in\R$   the sets $\{x+y\,|\,x\in
A,y\in B\}$ and $\{\lambda x\,|\,x\in A\}$ are denoted, as usually, by
$A+B$ and $\lambda A$.

Let $a$ be a Laurent polynomial in $n$ variables, $\Ge$  a positive
number and $\GG$  a face of the Newton polyhedron $\GD=\GD(a)$.

\begin{prop}\label{6.2.A}If in open set $U\subset K\R^n$ the truncation
$a^\GG$ is $\Ge$-sufficient for $a$, then it is $\Ge$-sufficient for $a$
in $qh_{\Cl DC_\GD^-(\GG)}(U)$.\footnote{Here (as above)
$\Cl$ denotes the closure.} \end{prop}

\begin{proof}Let $\Go\in\Cl DC_\GD^-(\GG)$ and $\Go w=\Gd$ for
$w\in\GG$. By \ref{6.1.E}, $\Ge$-sufficiency of truncation $a^\GG$
for $a$ in $qh_\Go(U)$ is equivalent to $\Ge$-sufficiency of
truncation $(a\circ qh_\Go)^\GG$ for $(a\circ qh_\Go)^\GG$ in $U$
or, equivalently, to $\Ge$-sufficiency of $\GG$-truncation of Laurent
polynomial $b=e^{-\Gd}a\circ qh_\Go$ in $U$. Since $$e^{-\Gd}a\circ
qh_\Go(x)=\sum_{w\in\GD}e^{-\Gd}a_wx^we^{\Go
w}=a^\GG(x)+\sum_{w\in\GD\sminus\GG}e^{\Go
w-\Gd}a_wx^w$$
and $\Go w-\Gd\le0$ when $w\in\GD\sminus\GG$ and
$\Go\in \Cl DC_\GD^-(\GG)$, it follows that $b$ satisfies the
conditions $\GD(b)=\GD$, $b^\GG=a^\GG$ and
$||b-b^\GG||\le||a-a^\GG||$. Therefore the truncation $b^\GG$ is
$\Ge$-sufficient for $b$ in $U$ and, hence, the truncation $a^\GG$
is $\Ge$-sufficient for $a$ in $qh_\Go(U)$. From this, by \ref{6.1.A},
the proposition follows.
\end{proof}

\begin{prop}\label{6.2.B}If the truncation $a^\GG$ is completely
nondegenerate and
$la V_{K\R^n}(a^\GG)$ has a tubular
\nei of the chosen type, then for any compact sets $C\subset K\R^n$ and
$\GO\subset DC_\GD^-(\GG)$ there exists $\Gd$ such that
in $qh_{\Gd\GO}(C)$ the truncation $a^\GG$ is
$\Ge$-sufficient for $a$.
\end{prop}

\begin{proof}For $\Go\in DC_\GD^-(\GG)$ denote by $\Go\GG$ a value
taken by the scalar product $\Go w$ for $w\in\GG$. Since
$$t^{-\Go\GG}a\circ
qh_{\Go,t}(x)=a^\GG(x)+\sum_{w\in\GD\sminus\GG}t^{\Go
w-\Go\GG}a_wx^w$$
for $\Go\in DC^-_\GD(\GG)$ (cf. the previous proof) and $\Go
w-\Go\GG<0$ when $w\in\GD\sminus\GG$ and $\Go\in
DC_\GD^-(\GG)$ it follows that the Laurent polynomial
$b_{\Go,t}=t^{-\Go\GG}a\circ qh_{\Go,t}$ with $\Go\in
DC_\GD^-(\GG)$ turns to $a^\GG$ as $\GG\to +\infty$. It is clear
that this convergence is uniform with respect to $\Go$ on a compact set
$\GO\subset DC_\GD^-(\GG)$. By \ref{6.1.D} it follows from this that
for a compact set $U\subset K\R^n$ there exists
$\eta$ such that for any $\Go\in\GO$ and $t\ge\eta$ the
truncation $b^\GG_{\Go,t}$ of
$b_{\Go,t}$ is $\Ge$-sufficient in $U$ for $b_{\Go,t}$. By
\ref{6.1.E}, the latter is equivalent to $\Ge$-sufficiency of
truncation $a^\GG$ for $a$ in $qh_{\Go,t}(U)$.

Thus if $U$ is the closure of a bounded \nei $W$ of a set $C$
then there exists $\eta$ such that for $\Go\in\GO$ and
$t\ge\eta$ the truncation $a^\GG$ is $\Ge$-sufficient for $a$ in
$qh_{\Go,t}(U)$. Therefore $a^\GG$ is the same in a smaller set
$qh_{\Go,t}(W)$ and, hence, (by \ref{6.1.A}) in the union
$\bigcup_{t\ge\eta,\Go\in\GO}qh_{\Go,t}(W)$ and, hence, in a smaller
set $\bigcup_{t=\eta,\Go\in\GO}(C)$. Putting
$\Gd=\ln\eta$ we obtain the required result.
\end{proof}

\begin{prop}\label{6.2.C}Let $\GG$ is a face of another face $\GS$ of
the polyhedron $\GD$. Let $\GO$ is a compact subset of the  cone
$DC^-_\GD(\GS)$. If $\GG$-truncation $a^\GG$ is $\Ge$-sufficient for
$a^\GS$ in a compact set $C$, then there exists a number $\Gd$
such that $a^\GG$ is $\Ge$-sufficient for $a$ in $qh_{\Gd\GO}(C)$.
\end{prop}

This proposition is proved similarly to \ref{6.2.B}, but with the
following difference: the reference to Theorem \ref{6.1.D} is replaced
by a reference to Theorem \ref{6.1.C}.\qed

\begin{prop}\label{6.2.D}Let $C\subset K\R^n$ be a compact set and let
$\GG$ be a face of $\GD$ such that
for any face $\GS$ of $\GD$ with $\dim\GS=\dim\GD-1$ having a face
$\GG$ the truncation $a^\GG$ is $\Ge$-sufficient for $a^\GS$ in $C$.
Then there exists a real number $\Gd$ such that the truncation $a^\GG$
is $\Ge$-sufficient for $a$ in $qh_{\Cl DC^-_\GD(\GG)\sminus \frak
N_\Gd DC^-_\GD(\GD)}(\Int C)$.  \end{prop}

\begin{proof}By \ref{6.2.C}, for any face $\GS$ of $\GD$ with
$\dim\GS=\dim\GD-1$ and $\GG\subset\p\GS$ there exists a
vector $\Go_\GS\in DC_\GD^-(\GS)$ such that the truncation
$a^\GG$ is $\Ge$-sufficient for $a$ in $qh_{\Go_\GS}(C)$,
and, hence, by \ref{6.2.A}, in $qh_{\Go_\GS+\Cl DC^-_\GD(\GG)}(\Int C)$.
Choose such $\Go_\GS$ for each $\GS$ with $\dim\GS=\dim\GD-1$ and
$\GG\subset\p\GS$. Obviously, the sets
$\Go_\GS+\Cl DC^-_\GD(\GG)$ cover the whole closure of the cone
$DC_\GD^-(\GG)$ besides some \nei of its top, i.e. the cone
$DC^-_\GD(\GD)$; in other words, there exists a number $\Gd$
such that
$\bigcup_\GS(\Go_\GS+\Cl DC_\GD^-(\GG))\supset
\Cl DC_\GD^-(\GG)\sminus \frak N_\Gd DC_\GD^-(\GD)$.
Hence, $a^\GG$ is $\Ge$-sufficient for $a$ in
$qh_{\Cl DC^-_\GD(\GG)\sminus \frak N_\Gd
DC^-_\GD(\GD)}(\Int C)\subset\bigcup_\GS
qh_{\Go_\GS+\Cl DC_\GD^-(\GG)}(\Int C)$.
\end{proof}

\subsection{The main Theorem on logarithmic asymptotes of
hypersurface}\label{s6.3}Let $\GD\subset\R^n$ be a convex closed
polyhedron and $\Gf:\mathcal{G}(\GD)\to\R$ be a positive function. Then for
$\GG\in\mathcal{G}(\GD)$ denote by $D_{\GD,\Gf}(\GG)$ the set $$\frak
N_{\Gf(\GG)}(DC_\GD^-(\GG))\sminus\bigcup_{\GS\in\mathcal{G}(\GD), \ 
\GG\in\mathcal{G}(\GS)} \frak N_{\Gf(\GS)}(DC_\GD^-(\GS)).$$

It is clear that the sets $D_{\GD,\Gf}(\GG)$ with $\GG\in\mathcal{G}(\GD)$ cover
$\R^n$. Among these sets only sets corresponding to faces of the
same dimension can intersect each other. In some cases (for example,
if $\Gf^\GG$ grows fast enough when $\dim\GG$ grows) they do not
intersect  and then $\{D_{\GD,\Gf}(\GG)\}_{\GG\in\mathcal{G}(\GD)}$ is a
partition of $\R^n$.

Let $a$ be a Laurent polynomial over $K$ in $n$ variables and $\Ge$ be a
positive number. A function $\Gf:\mathcal{G}(\GD(a))\to\R^n$ is said to be
{\it describing domains of $\Ge$-sufficiency for }
$a$ (with respect to the chosen class of tubular neighborhoods) if for
any proper face $\GG\in\mathcal{G}(\GD(a))$, for which truncation $a^\GG$ is
completely non-degenerate and the hypersurface $la(V_{K\R^n}(a^\GG))$
has a tubular \nei of the chosen class, the truncation $a^\GG$ is
$\Ge$-sufficient for $a$ in some \nei of $l^{-1}(DC_{\GD(a),\Gf}(\GG))$.

\begin{prop}\label{6.3.A}For any Laurent polynomial $a$ over $K$ in
$n$ variables and $\Ge>0$ there exists a function
$\mathcal{G}(\GD(a))\to\R$ describing domains of $\Ge$-sufficiency for $a$
with respect to the chosen class of tubular neighborhoods.
\end{prop}

In particular, if $a$ is peripherally nondegenerate Laurent polynomial
over $K$ in $n$ variables and $\dim\GD(a)=n$ then for any $\Ge>0$ there
exists a compact set $C\subset K\R^n$ such that $K\R^n\sminus C$ is
covered by regions in which truncations of
$a^{\p\GD(a)}$ are $\Ge$-sufficient for $a$ with respect
to class of normal tubular neighborhoods. In other words, under these
conditions behavior of $V_{K\R^n}(a)$ outside $C$ is defined by
monomials of $a^{\p\GD(a)}$.

\subsection{Proof of Theorem \ref{6.3.A}}\label{s6.4}Theorem \ref{6.3.A} is
proved by induction on dimension of polyhedron $\GD(a)$.

If $\dim\GD(a)=0$ then $a$ is monomial and $V_{K\R^n}(a)=\varnothing$.
Thus for any $\Ge>0$ any function $\Gf:\mathcal{G}(\GD(a))\to\R$ describes domains of $\Ge$-sufficiency for $a$.

Induction step follows obviously from the following Theorem.

\begin{prop}\label{6.4.A}Let $a$ be a Laurent polynomial over $K$ in $n$
variables, $\GD$ be its Newton polyhedron, $\Ge$ a positive number.
If for a function $\Gf:\mathcal{G}(\GD)\sminus\{\GD\}\to\R$ and
any proper face $\GG$ of $\GD$ the restriction $\Gf|_{\mathcal{G}(\GG)}$
describes domains of $\Ge$-sufficiency for $a^\GG$, then $\Gf$
can be extended to a function $\bar\Gf:\mathcal{G}(\GD)\to\R$
describing regions of $\Ge$-sufficiency for $a$.
\end{prop}

\begin{proof}It is sufficient to prove that for any face
$\GG\in\mathcal{G}(\GD)\sminus\{\GD\}$, for which the truncation
$a^\GG$ is completely
nondegenerate and hypersurface $V_{K\R^n}(a^\GG)$ has a
tubular \nei of the chosen class, there exists an extension
$\Gf_\GG$ of $\Gf$ such that truncation $a^\GG$ is $\Ge$-sufficient
for $a$ in a \nei of $l^{-1}(D_{\GD,\Gf_\GG}(\GG))$, i.e. to
prove that for any face $\GG\ne\GD$ there exists a number
$\Gf_\GG(\GD)$ such that the truncation $a^\GG$ is $\Ge$-sufficient for
$a$ in some \nei of
$$l^{-1}(\frak N_{\Gf(\GG)}(DC^-_\GD(\GG))\sminus[\frak
N_{\Gf_\GG(\GD)}(DC^-_\GD(\GD))\bigcup_{\GS\in\mathcal{G}(\GD)\sminus\{\GD\}, \ 
\GG\in\mathcal{G}(\GS)}\frak
N_{\Gf(\GS)}(DC^-_\GD(\GS))].$$
Indeed, putting
$$\bar\Gf(\GD)=\max_{\GG\in\mathcal{G}(\GD)\sminus\{\GD\}}\Gf_\GG(\GD)$$
we obtain a required
extension of $\Gf$.

First, consider the case of a face $\GG$ with
$dim\GG=dim\GD-1$. Apply proposition \ref{6.2.B} to $C=l^{-1}(\Cl\frak
N_{\Gf(\GG)+1}(0)$ and any one-point set $\GO\subset DC^-_\GD(\GG)$.
It implies that $a^\GG$ is $\Ge$-sufficient for $a$ in $qh_\Go(C)=
l^{-1}(Cl\frak N_{\Gf(\GG)+1}(\Go))$ for some
$\Go\in DC^-_\GD(\GG)$. Now  apply proposition \ref{6.2.A} to $U=
l^{-1}(\frak N_{\Gf(\GG)+1}(\Go))$. It gives that $a^\GG$ is
$\Ge$-sufficient for $a$ in $qh_{DC^-_\GD(\GG)}(
l^{-1}(\frak N_{\Gf(\GG)+1}(\Go))=l^{-1}(\frak N_{\Gf(\GG)+1}
(\Go+DC^-_\GD(\GG)))$ and, hence, in the smaller set
$l^{-1}(\frak N_{\Gf(\GG)+1}(DC^-_\GD(\GG)))\sminus\frak
N_{|\Go|}(DC^-_\GD(\GD))$. It is remained to put
$\Gf_\GG(\GD)=|\Go|+1$.

Now consider the case of face $\GG$ with $\dim\GG<\dim\GD-1$. Denote by
$E$ the set
$$\frak N_{\Gf(\GG)}(DC^-_\GD(\GG))\sminus\bigcup_{\GS\in\mathcal{G}(\GD)\sminus\{\GD\}, \ \GG\in\mathcal{G}(\GS)}\frak
N_{\Gf(\GS)}(DC^-_\GD(\GS)).$$
It is clear that there exists a ball $B\subset\R^n$ with center at
$0$ such that $E=(E\cap B)+\Cl DC^-_\GD(\GG)$. Denote
the radius of this ball by $\Gb$.

If $\GS\in\mathcal{G}(\GD)$ is a face of dimension $\dim\GD-1$ with
$\p\GS\supset\GG$ then, by the hypothesis, the truncation
$a^\GG$ is $\Ge$-sufficient for $a^\GS$ in some \nei of
$$l^{-1}(\frak
N_{\Gf(\GG)}(DC^-_\GS(\GG))\sminus\bigcup_{\Theta\in\mathcal{G}(\GS), \ \GG\in\mathcal{G}(\Theta)}\frak
N_{\Gf(\Theta)}(DC^-_\GS(\Theta))$$
and, hence, in \nei of a smaller set
$$l^{-1}(\frak
N_{\Gf(\GG)}(DC^-_\GD(\GG))\sminus\bigcup_{\Theta\in\mathcal{G}(\GS), \ \GG\in\mathcal{G}(\Theta)}\frak
N_{\Gf(\Theta)}(DC^-_\GD(\Theta)).$$
Therefore for any face $\GS$ with $\dim\GS=\dim\GD-1$ and
$\GG\subset\p\GS$ the truncation $a^\GG$ is $\Ge$-sufficient for
$a^\GS$ in some \nei of $l^{-1}(E)$. Denote by $C$ a compact
\nei of $l^{-1}(E\cap B)$ contained in this neighborhood. Applying
proposition \ref{6.2.A}, one obtains that $a^\GG$ is $\Ge$-sufficient
for $a$ in the set $$qh_{\Cl DC^-_\GD(\GG)\sminus\frak N_\Gd
DC^-_\GD(\GD)}(\Int C)=l^{-1}(\Int
l(C)+\Cl DC^-_\GD(\GG)\sminus\frak N_\Gd
DC^-_\GD(\GD))).$$
It is remained to put $\Gf_\GG(\GD)=\Gd+\beta$
\end{proof}

\subsection{Modification of Theorem 6.3.A}\label{s6.5}Below in Section
\ref{s6.8} it will be more convenient to use not Theorem \ref{6.3.A}
but the following its modification, whose formulation is more cumbrous,
and whose proof is obtained by an obvious modification of deduction of
\ref{6.3.A} from \ref{6.4.A}.

\begin{prop}\label{6.5.A}For any Laurent polynomial $a$ over $K$ in $n$
variables and any $\Ge>0$ and $c>1$ there exists a function 
$\Gf:\mathcal{G}(\GD(a))\to\R$ such that for any proper face $\GG\in\mathcal{G}(\GD(a))$, for which $a^\GG$ is completely nondegenerate and
$la(V_{K\R^n}(a^\GG))$ has a tubular \nei
from the chosen class, the truncation $a^\GG$ is $\Ge$-sufficient for
$a$ in some \nei of
$$l^{-1}(\frak
N_{c\Gf(\GG)}(DC^-_\GD(\GG))\sminus\bigcup_{\GS\in\mathcal{G}(\GD), \ \GG\in\mathcal{G}(\GS)}\frak
N_{\Gf(\GS)}(DC^-_\GD(\GS)).$$
\qed\end{prop}

\subsection{Charts of L-polynomials}\label{s6.6}Let $a$ be a
peripherally nondegenerate Laurent polynomial over $K$ in $n$
variables, $\GD$ be its Newton polyhedron. Let $V$ be a vector subspace
of $\R^n$ corresponding to the smallest affine subspace containing
$\GD$ (i.e.  $V=C_\GD(\GD))$. Let $\Gf:\mathcal{G}(\GD)\to\R$ be the
function, existing by \ref{6.3.A}, describing for some $\Ge$ regions of
$\Ge$-sufficiency for $a$ with respect to class of normal tubular
neighborhoods.

The pair $(\GD\times U_K^n,$ $\upsilon)$ consisting of the product
$\GD\times U_K^n$ and its subset $\upsilon$ is a $K$-{\it chart} of a
Laurent polynomial $a$ if:
\begin{enumerate}
\item\label{6.6.A}there exists a homeomorphism
$h:(\Cl D_{\GD,\Gf}(\GD)\cap V)\times U_K^n\to\GD\times U_K^n$
such that $h((Cl\phantom{a}D_{\GD,\Gf}(\GD)\cap V)\times
y)=\GD\times y$ for $y\in U_K^n$,
$$\upsilon=h(la V_{K\Bbb
R^n}(a)\cap(\Cl D_{\GD,\Gf}(\GD)\cap V)\times U_K^n$$
and
for each face $\GG$ of $\GD$ the set $h((\Cl D_{\GD,\Gf}(\GD)
\cap D_{\GD,\Gf}(\GG)\cap V)\times U_K^n)$ lies in the product of the
star $\bigcup_{\GG\in\mathcal{G}(\GS)\\\GS\in\mathcal{G}(\GD)\sminus\{\GD\}}\GS$ of $\GG$ to $U_K^n$;

\item\label{6.6.B} for any vector $\Go\in\R^n$, which is
orthogonal to $V$ and, in the case of $K=\R$, is integer, the set
$\upsilon$ is invariant under transformation $\GD\times U_K^n\to
\GD\times U_K^n$ defined by formula $(x,$ $(y_1,$\dots, $y_n))\mapsto
(x,$ $(e^{\pi i\Go_1}y_1,$ \dots, $e^{\pi i\Go_n}y_n))$;

\item\label{6.6.C}for each face $\GG$ of $\GD$ the pair $(\GG\times
U_K^n,$ $\upsilon\cap(\GG\times U_K^n))$ is a $K$-chart of Laurent
polynomial $a^\GG$.
\end{enumerate}

The definition of the chart of a Laurent polynomial, which, as I
believe, is clearer than the description given here, but based on the
notion of toric completion of $K\R^n$, is given above in Section
\ref{s3.3}.  I restrict myself to the following commentary of
conditions \ref{6.6.A} -- \ref{6.6.C}.

The set $(\Cl D_{\GD,\Gf}(\GD)\cap V)\times U_K^n$ contains, by
\ref{6.3.A}, a deformation retract of $la V_{K\R^n}(a)$. Thus,
condition \ref{6.6.A} means that $\upsilon$ is homeomorphic to a deformation
retract of $V_{K\R^n}(a)$. The position of $\upsilon$
in $\GD\times U_K^n$ contains, by \ref{6.6.A} and \ref{6.6.C}, a
complete topological information about behavior of this hypersurface
outside some compact set. The meaning of \ref{6.6.B} is in that $\upsilon$
has the same symmetries as, according to \ref{2.1.C}, $V_{K\R^n}(a)$
has.

\subsection{Structure of $V_{K\R^n}(b_t)$ with small $t$}\label{s6.7}
Denote by $i_t$ the embedding $K\R^n\to K\R^{n+1}$ defined
by  $i_t(x_1,\dots, x_n)=(x_1, \dots,x_n,\,t)$. Obviously,
$$V_{K\R^n}(b_t)=i^{-1}_t V_{K\R^{n+1}}(b).$$
This allows to take advantage of results of the previous
Section for study of $V_{K\R^n}(b_t)$ as $t\to0$.
For sufficiently small $t$   the  image of embedding $i_t$ is covered by
regions of $\Ge$-sufficiency of truncation $b^{\tilde\GG}$, where
$\tilde\GG$ runs over the set of faces of graph of  $\nu$, and
therefore the hypersurface $V_{K\R^n}(b_t)$ turns
to be composed of pieces obtained from
$V_{K\R^n}(a_i)$ by appropriate quasi-homotheties.

I preface the formulation describing in detail the behavior of
$V_{K\R^n}(b_t)$ with several notations.

Denote  the Newton polyhedron $\GD(b)$ of Laurent
polynomial $b$ by $\tilde\GD$. It is clear that $\tilde\GD$ is
the convex hull of the
graph of $\nu$. Denote by $\mathcal{G}$ the union
$\bigcup_{i=1}^s\mathcal{G}(\GD_i)$. For $\GG\in\mathcal{G}$ denote by $\tilde\GG$
the graph $\nu|_\GG$. It is clear that $\tilde\GG\in\mathcal{G}(\tilde\GD)$
and hence an injection $\GG\to\tilde\GG:\mathcal{G}\to\mathcal{G}(\tilde\GD)$ is
defined.

For $t>0$ denote by $j_t$ the embedding $\R^n\to\R^{n+1}$
defined by the formula $j_t(x_1, \dots,x_n)=(x_1, \dots,x_n,\,\ln t)$.
Let $\psi:\mathcal{G}\to\R$ be a positive function, $t$ be a number from
interval $(0,1)$. For $\GG\in\mathcal G$ denote by $\mathcal{E}_{t,\psi}(\GG)$ the
following subset of $\R^n$:
$$\frak
N_{\psi(\GG)}j_t^{-1}(DC^-_{\tilde\GD}(\tilde\GG))\sminus
\bigcup_{\GS\in\mathcal{G}, \ \GG\in\mathcal{G}(\GS)}\frak
N_{\Gf(\GS)}j_t^{-1}(DC^-_{\tilde\GD}(\tilde\GS)).$$

\begin{prop}\label{6.7.A}If Laurent polynomials $a_1, \dots, a_s$ are
completely non-degenerate then for any $\Ge>0$ there exist $t_0\in(0,1)$
and function $\psi:\mathcal{G}\to\R$ such that for any $t\in(0,t_0]$ and
any face $\GG\in\mathcal G$ truncation $b_t^\GG$ is
$\Ge$-sufficient for $b_t$ with respect to the class of normal
tubular \neis in some \nei of $l^{-1}(\mathcal{E}_{t,\psi}(\GG))$.
\end{prop}

Denote  the gradient of restriction of $\nu$ on
$\GG\in\mathcal G$ by $\nabla(\GG)$. The truncation $b_t^\GG$, obviously,
equals $a^\GG\circ qh_{\nabla(\GG),t}$. In particular,
$b_t^{\GD_i}=a_i\circ qh_{\nabla(\GD_i),t}$ and, hence,
$$V_{K\R^n}(b_t^{\GD_i})=qh_{\nabla(\GD_i),t^{-1}}(V_{K\Bbb
R^n}(a_i)).$$

In the domain, where $b^\GG_t$ is $\Ge$-sufficient for $b_t$, the
hypersurfaces $la V_{K\R^n}(b_t)$ and $la V_{K\R^n}(b_t^{\GD_i})$
with $\GD_i\supset\GG$ lie in the same normal tubular
$\Ge$-\nei of $la V_{K\R^n}(b_t^\GG)$ and, hence, are isotopic by an
isotopy moving points at most on $2\Ge$. Thus, according to \ref{6.7.A},
for $t\le t_0$ to the space $K\R^n$ is covered by regions in which
$V_{K\R^n}(b_t)$ is approximated by
$qh_{\nabla(\GD_i),t^{-1}}(V_{K\R^n}(a_i))$.

\subsection{Proof of Theorem \ref{6.7.A}}\label{s6.8}Put
$c=\max\{\sqrt{1+\nabla(\GD_i)^2},$ $i=1,\dots, s\}$. Apply Theorem
\ref{6.5.A} to the Laurent polynomial $b$ and numbers $\Ge$ and $c$,
considering as the class of chosen tubular \neis in $\R^{n+1}\times
U_K^{n+1}$ tubular neighborhoods, whose fibers lie in the fibers
$\R^n\times t\times U_K^n\times s$ of $\R^{n+1}\times U_K^{n+1}$ and
consist of segments of geodesics which are orthogonal to intersections
of submanifold with $\R^n\times t\times U_K^n\times s$. (Intersection
of such a tubular \nei of $M\subset\R^{n+1}\times U_K^{n+1}$ with the
fiber $\R^n\times t\times U_K^n\times s$ is a normal tubular \nei of
$M\cap(\R^n\times t\times U_K^n\times s)$ in $\R^n\times
t\times U_K^n\times s$.) Applying  Theorem \ref{6.5.A} one
obtains a function $\Gf:\mathcal{G}(\tilde\GD)\to\R$. Denote by $\psi$
the function $\mathcal{G}\to\R$ which is the composition of embedding
$\GG\mapsto\tilde\GG:\mathcal{G}\to\mathcal{G}(\tilde\GD)$ (see Section
\ref{s6.7}) and the function $\dfrac1c\Gf:\mathcal{G}(\tilde\GD)\to\R$. This
function has the required property. Indeed, as it is easy to see,
for $0<t<e^{-\Gf(\tilde\GD)}$
$\mathcal{E}_{t,\psi}(\GG)$ is contained,  in
$$j_t^{-1}(\frak
N_{c\Gf(\tilde\GG)}(DC^-_{\tilde\GD}(\tilde\GG))\sminus
\bigcup_{\tilde\GS\in\mathcal{G}(\tilde\GD), \ \tilde\GG\in\mathcal{G}(\tilde\GS)}\frak
N_{\Gf(\tilde\GS)}(DC^-_{\tilde\GD}(\tilde\GS))),$$
and thus from $\Ge$-sufficiency of $b^{\tilde\GG}$ for $b$ in  some \nei
of
$$l^{-1}(\frak
N_{c\Gf(\tilde\GG)}(DC^-_{\tilde\GD}(\tilde\GG))\sminus
\bigcup_{\tilde\GS\in\mathcal{G}(\tilde\GD), \ \tilde\GG\in\mathcal{G}(\tilde\GS)}\frak
N_{\Gf(\tilde\GS)}(DC^-_{\tilde\GD}(\tilde\GS))),$$
with respect  to the chosen class of tubular \neis in $\R^{n+1}\times
U_K^{n+1}$ if follows that for $0<t<e^{-\Gf(\tilde\GD)}$ the truncation
$b_t^\GG$ is $\Ge$-sufficient  for $b_t$ in some \nei of
$l^{-1}(\mathcal{E}_{t,\psi}(\GG))$ with respect to the class of normal
tubular neighborhoods.  \qed

\subsection{An alternative proof of Theorem \ref{4.3.A}}\label{s6.9}Let
$V$ be a vector subspace of $\R^n$ corresponding to the minimal affine
subspace containing $\GD$. It is divided for each $t\in(0,1)$ onto
the sets $V\cap j_t^{-1}(DC^-_{\tilde\GD}(\tilde\GG))$ with $\GG\in\mathcal{G}$. 
Let us construct cells $\GG_t$ in $V$ which are dual to the sets of
this partition (barycentric stars). For this mark a point in each $V\cap
j_t^{-1}(DC^-_{\tilde\GD}(\tilde\GG))$:
$$b_{t,\GG}\in
V\cap j_t^{-1}(DC^-_{\tilde\GD}(\tilde\GG)).$$
Then for $\GG$ with
$\dim\GG=0$ put $\GG_t=b_{t,\GG}$ and construct  the others $\GG_t$
inductively on dimension $\dim\GG$: if $\GG_t$ for $\GG$ with $\dim\GG<r$
have been constructed then for $\GG$ with $\dim\GG=r$ the cell $\GG_t$
is the (open) cone on $\bigcup_{\GS\in\mathcal{G}(\GG)\sminus\{\GG\}}\GS_t$
with the vertex $b_{t,\GG}$. (This is the usual construction of dual
partition turning in the case of triangulation to partition onto
barycentric stars of simplices.)

By Theorem \ref{6.7.A} there exist $t'_0\in(0,1)$ and function $\psi:\mathcal{G}\to\R$ such that for any $t\in(0,t'_0]$ and any face $\GG\in\mathcal G$
the truncation $b_t^\GG$ is $\Ge$-sufficient for $b_t$ in some \nei of
$l^{-1}(\mathcal{E}_{t,\psi}(\GG))$. Since cells $\GG_t$ grow unboundedly
when $t$ runs to zero (if $\dim\GG\ne0$) it follows that there exists
$t_0\in(0,t'_0]$ such that for $t\in(0,t_0]$ for each face $\GG\in\mathcal
G$ the set $\frak N_{\psi(\GG)}j^{-1}_t(DC^-_{\tilde\GD}(\GG))$, and
together with it the set $\mathcal{E}_{t,\psi}(\GG)$, lie in the star of
the cell $\GG_t$, i.e. in $\bigcup_{\GG\in\mathcal{G}(\GS)}\GS_t$. Let us
show that for such $t_0$ the conclusion of Theorem \ref{4.3.A} takes place.

Indeed, it follows from \ref{6.7.A} that there exists a homeomorphism
$h:\GG_t\times U_K^n\to\GG\times U_K^n$ with $h(\GG_t\times y)=\GG\times y$
for $y\in U_K^n$ such that $(\GG\times U_K^n,$ $h(la(V_{K\Bbb
R^n}(b_t))\cap\GG_t\times U_K^n))$ is $K$-chart of Laurent polynomial
$a^\GG$. Therefore the pair
$$(\cup_{\GG\in\mathcal{G}}\GG_t\times U_K^n,\,
laV_{K\R^n}(b_t)\cap(\cup_{\GG\in\mathcal{G}}\GG_t\times U_K^n))$$
is
obtained in result of patchworking $K$-charts of Laurent polynomials
$a_1,$ \dots, $a_s$. The function $\Gf:\mathcal{G}(\GD)\to\R$, existing
by Theorem \ref{6.3.A} applied to $b_t$, can be chosen, as it follows from
\ref{6.4.A}, in such a way that it should majorate any given in advance
function $\mathcal{G}(\GD)\to\R$. Choose $\Gf$ in such a way that
$D_{\Gf,\GD}(\GD)\supset\cup_{\GG\in\mathcal{G}}\GG_t$ and
$D_{\Gf,\GD}(\GS)\cap\p D_{\Gf,\GD}(\GD)\supset\mathcal{E}_{t,\psi}(\GS)
\cap\p D_{\Gf,\GD}(\GD)$ for $\GS\in\mathcal{G}(\GD)\sminus\{\GD\}$. As it follows from \ref{6.7.A}, there exists a
homeomorphism
\begin{multline}(\bigcup_{\GG\in\mathcal{G}}\GG_t\times U_K^n,\,la
V_{K\R^n}(b_t)\cap(\bigcup_{\GG\in\mathcal{G}}\GG_t\times
U_K^n))\to\\(D_{\Gf,\GD}(\GD)\times U_K^n,\,la V_{K\R^n}(b_t)\cap(
D_{\Gf,\GD}(\GD)\times U_K^n))\end{multline}
turning $\mathcal{E}_{t,\psi}(\GS)\cap\p(\bigcup_{\GG\in\mathcal{G}}\GG_t\times U_K^n)$ to $\mathcal{E}_{t,\psi}(\GS)\cap\p
D_{\Gf,\GD}(\GD)$ for $\GS\in\mathcal{G}(\GD)\sminus\{\GD\}$. Therefore
$K$-chart of Laurent polynomial $b_t$ is obtained by patchworking
$K$-charts of Laurent polynomials $a_1, \dots, a_s$.
\qed


\bibliographystyle{amsalpha}
\bibliography{crav}
\end{document}